\documentclass[11pt]{article}
\usepackage{amsmath}
\usepackage{amssymb}
\usepackage{amsthm}
\usepackage{amsfonts}
\usepackage{bm}
\usepackage{braket}
\usepackage[truedimen,margin=25truemm]{geometry}
\usepackage[dvipdfmx]{hyperref}

\allowdisplaybreaks

\theoremstyle{plain}
\newtheorem{th.}{Theorem}[section]
\newtheorem{prop.}[th.]{Proposition}
\newtheorem{lem.}[th.]{Lemma}
\newtheorem{cor.}[th.]{Corollary}
\newtheorem{def.}[th.]{Definition}
\newtheorem{rmk.}[th.]{Remark}
\newtheorem{conj.}[th.]{Conjecture}

\newcommand{\BA}{\mathbb{A}}

\newcommand{\BC}{\mathbb{C}}

\newcommand{\BE}{\mathbb{E}}

\newcommand{\BP}{\mathbb{P}}

\newcommand{\BR}{\mathbb{R}}
\newcommand{\BS}{\mathbb{S}}
\newcommand{\BT}{\mathbb{T}}

\newcommand{\CA}{\mathcal{A}}
\newcommand{\CB}{\mathcal{B}}

\newcommand{\CF}{\mathcal{F}}

\newcommand{\CN}{\mathcal{N}}

\newcommand{\CS}{\mathcal{S}}

\newcommand{\ga}{\gamma}
\newcommand{\de}{\delta}
\newcommand{\ep}{\varepsilon}
\newcommand{\ps}{\psi}
\newcommand{\ph}{\varphi}
\newcommand{\ta}{\tau}

\renewcommand{\th}{\theta}
\newcommand{\si}{\sigma}
\newcommand{\om}{\omega}

\newcommand{\et}{\eta}

\newcommand{\rh}{\rho}
\newcommand{\ch}{\chi}

\newcommand{\De}{\Delta}
\newcommand{\Ph}{\Phi}
\newcommand{\Ps}{\Psi}

\newcommand{\Om}{\Omega}

\newcommand{\pl}{\partial}
\newcommand{\wt}{\widetilde}
\newcommand{\wh}{\widehat}

\newcommand{\Dk}[1]{\left[#1\right]}
\newcommand{\K}[1]{\left(#1\right)}
\newcommand{\Br}[1]{\langle #1 \rangle}
\newcommand{\Abs}[1]{\left|#1\right|}
\newcommand{\No}[1]{\left\| #1 \right\|}
\newcommand{\AL}[1]{\begin{align} #1 \end{align}}
\newcommand{\ALN}[1]{\begin{align*} #1 \end{align*}}

\newcommand{\I}{\infty}

\newcommand{\sd}{\langle \nabla \rangle}

\newcommand{\bx}{\langle x \rangle}

\newcommand{\lxr}{\langle \xi \rangle}

\renewcommand{\H}{{\dot{H}^{1/4}_x \cap \dot{H}^{1/4,\I}_x}}
\newcommand{\Ltx}{{L^2_tL^2_x}}

\newcommand{\na}{\nabla}

\renewcommand{\Re}{\mathrm{Re}}

\newcommand{\tw}{\frac{1}{2}}
\newcommand{\ft}{\frac{1}{4}}
\newcommand{\ov}{\overline}
\newcommand{\dH}{\dot{H}}

\title{Asymptotic stability of a wide class of steady states for the Hartree equation for random fields}
\author{Sonae HADAMA
\thanks{Kyoto University.
E-mail address: \texttt{hadama@kurims.kyoto-u.ac.jp}}
}
\date{}

\begin{document}
\maketitle
\vspace{-5mm}
\begin{abstract}We study the Hartree equation describing the time evolution of the wave functions of infinitely many fermions interacting with each other.
	The Hartree equation can be formulated in terms of
	random fields.
	This formulation was introduced by de Suzzoni in \cite{2015S}.
	It has infinitely many steady states, and \cite{2020CS, 2022CS} have studied its asymptotic stability.
	However, they required somewhat strong conditions for steady states, and the stability of Fermi gas at zero temperature, one of the most important steady states from the physics point of view, was left open.
	In this paper, we prove the stability of steady states in a wide class, which includes Fermi gas at zero temperature in $d$-dimensional space when $d \ge 3$.
	
	This paper is the revised version.
	There are two major changes.
	The first one is that the assumptions for potential $w$ in Theorem 1.3 are weakened.
	The second one is that we removed Theorem 1.4 in the previous version because the author found some gaps which seemed difficult to fix.
\end{abstract}

\noindent {\bf 2020 MSC} Primary, 35Q40. Secondary, 35B35,  35B40.

\noindent {\bf Keywords---}Hartree equation, Asymptotic stability, Scattering, Fermi gas at zero temperature, Random fields, Density functional theory.

\tableofcontents

\section{Introduction}
In this paper, we consider the Hartree equation introduced in \cite{2015S}:
\begin{align}\label{NLH}\tag{NLH} 
		i\pl_t X = - \De X + (w \ast \BE[|X|^2])X,
\end{align}
where $X:\BR_t \times \BR^d_x \times \Om\to \BC$ is a random field defined over a probability space $(\Om, \CA, \BP)$ with expectation $\BE$. 
Here, $\ast$ is the space convolution and $w$ is a finite signed Borel measure on $\BR^d$.

We briefly explain the background of \eqref{NLH} for convenience according to \cite{2015S, 2015LS}; 
for more details, see \cite{2015S, 2020CS, 2015LS} and \cite{2022CS}.
The time evolution of wave functions of $N$ fermions in $d$-dimensional space is described by the linear Schr\"{o}dinger equation on $\BR^{Nd}$.
As one of the approximations to this equation, we have the following Hartree equation:
\begin{equation} \label{originalH}
	i \pl_t u_n(t,x) = \left(- \Delta_x +  w \ast \left( \sum_{m=1}^N |u_m(t,x)|^2 \right) \right) u_n(t,x),
	\quad (t,x) \in \BR \times \BR^d, \quad (n=1, \dots N).
\end{equation}
In this paper, we would like to deal with the case $N = \I$, that is, the Hartree equation for an infinite system of fermions.
However, it is difficult to study the case $N = \I$ in this form.
Thus, we rewrite (\ref{originalH}).
First we assume that $N < \I$.
Let $(u_n)_{n=1}^N \subset C(\BR, L^2_x)$ be a solution to \eqref{originalH}.
Let $(g_n)_{n=1}^N$ be an orthonormal family of $L^2_\om = L^2(\Om)$.
Define $X:\BR_t \times \BR^d_x \times \Om \to \BC$ by
\begin{align*}
	X(t,x,\om) := \sum_{n=1}^N u_n(t,x) g_n(\om).
\end{align*}
Then we can easily see that this $X$ satisfy \eqref{NLH}.
Note that \eqref{NLH} does not explicitly include the number of particles $N$.
Therefore, this equation is more suitable for our analysis, which deals with the case when $N=\I$.

It is known that \eqref{NLH} has infinitely many steady states.
More precisely, we have the following easy, but very important proposition.
\begin{prop.}[\cite{2015S}Proposition 1.6; See \cite{2020CS, 2022CS}]\label{SW}
 Let $d \ge 1$.
 Let $w$ be a finite signed Borel measure on $\BR^d$. 
 Let $f \in L^2_\xi$ and $m:= \hat{w}(0) \|f\|_{L^2_\xi}^2$.
 Define $Y_f$ by
 \begin{align*}
 	Y_f(t,x,\om):= \int_{\BR^d} e^{ix\xi - it(m+|\xi|^2)}f(\xi) dW(\xi),
 \end{align*}
 where $\int_{\BR^d} dW(\xi)$ is the Wiener integral.
 Then $Y_f$ is a solution to \eqref{NLH}.
 Moreover, the distribution of $Y_f$ does not depends on $(t,x)$,
 that is, $Y_f(t,x,\om) \sim \BC \CN(0,\|f\|_{L^2_\xi}^2)$,
 where $X \sim \BC \CN(0,\si^2)$ means the random variable $X$ follows the complex centered normal distribution with variance $\si^2$.
\end{prop.}

Proposition \ref{SW} implies that \eqref{NLH} has infinitely many steady states.
Furthermore, some $f$ are important from the physics point of view.
For example,
\begin{itemize}
	\item Fermi gas at zero temperature and chemical potential $\mu >0$:
	\AL{ \label{P1}
		|f(\xi)|^2 = \ch_{\{|\xi|^2 \leq \mu\}}(\xi);
	}
	\item Fermi gas at positive temperature $T>0$ and chemical potential $\mu \in \BR$:
	\AL{ \label{P2}
		|f(\xi)|^2 = \frac{1}{e^{(|\xi|^2-\mu)/T}+1};
	}
	\item Bose gas at positive temperature $T>0$ and chemical potential $\mu < 0$:
	\AL{ \label{P3}
		|f(\xi)|^2 = \frac{1}{e^{(|\xi|^2-\mu)/T}-1};
	}
	\item Boltzmann gas at positive temperature $T>0$ and chemical potential $\mu \in \BR$:
	\AL{ \label{P4}
		|f(\xi)|^2 = e^{-(|\xi|^2-\mu)/T};
	}
\end{itemize}
where $\ch_A$ is the indicator function of $A$.

In this paper, we consider the stability of $Y = Y_f$.
Before stating our main results, we discuss known results and related studies.
The first study of \eqref{NLH} is \cite{2015S}.
\cite{2015S} studied the Cauchy and scattering problems not only on $\BR^d$ but also on the torus $\BT^d$ and the sphere $\BS^d$.
\cite{2015S} proved analogues with basic results for the Hartree equation for one-particle.
The first study of the asymptotic stability of $Y$ is \cite{2020CS}, which proved the case when $d \ge 4$.
\cite{2022CS} extended the previous stability result to $d = 2,3$.
\cite{2022CS} succeeded to deal with potentials $w$ in a wide class that includes the delta measure.
However, \cite{2020CS, 2022CS} needed somewhat strong conditions for $f$,
in particular, $f$ has to be smooth.
Thus, their results do not include Fermi gas at zero temperature, that is, $|f(\xi)|^2 = \ch_{\{|\xi|\le \mu \}}$ with $\mu > 0$.

\eqref{NLH} is closely related to the Hartree equation for density matrix:
\begin{align} \label{gamma}
	i\pl_t \ga = [-\De+w \ast \rh_\ga , \ga], \quad \ga:\BR \to \CB(L^2_x),
\end{align}
where $\CB(L^2_x)$ is the set of all bounded linear operators on $L^2_x$, 
$[\cdot, \cdot]$ is the commutator,
$\rh_A(x) : \BR^d_x \to \BC$ is the density function of $A$, 
that is, $\rh_A(x):=k(x,x)$ with $k(x,y)$ is the integral kernel of $A$.
It is known that Fourier multiplier $g(-i\na) = \CF^{-1} g \CF \in \CB(L^2_x)$ is a stationary solution to \eqref{gamma}
for any $g \in L^1_\xi \cap L^\I_\xi$.
$|f|^2(-i\na)$ for \eqref{gamma} is a counterpart of $Y_f$ for \eqref{NLH}.
For other detailed parallel between \eqref{NLH} and \eqref{gamma}, see \cite{2015S, 2020CS}.

The Cauchy and scattering problems for \eqref{gamma} itself have been studied by \cite{1974BDF, 1976BDF, 1976C, 1992Z} and more recently, by \cite{2021PS}.
(More precisely, they considered more general nonlinear terms.)
However, the first study of stationary solution $g(-i\na)$ for \eqref{gamma} is the seminal works \cite{2015LS} and \cite{2014LS}.
\cite{2015LS} considered the Cauchy problem,
and \cite{2014LS} proved the asymptotic stability of $g(-i\na)$ when $d=2$.
One of the most important tools for their analysis is the Strichartz estimates generalized to orthonormal systems,
which is proven \cite{2014FLLS}.
For this topic, see also \cite{2017FS, 2019BHLNS}.
In relation to \cite{2014LS, 2015LS}, \cite{2020LS} proved the semi-classical limit of the Hartree equation is the Vlasov eqation:
\begin{align*}
	\pl_t W + 2v \cdot \na_x W - \na_x(w \ast \rh_W) \cdot \na_v W = 0,
\end{align*}
where $\rh_W=(2\pi)^{-d} \int  W(v,x) dv$.
The initial data was taken from near the stationary solution given by the form $g(-i\na)$.
And they get the global well-posed result for Vlasov equation as a by-product.

\cite{2017CHP} extended the global well-posedness results of \cite{2015LS}.
\cite{2018CHP} proved the asymptotic stability of $g(-i\na)$ when $d\ge3$.
However, \cite{2014LS, 2018CHP} needed somewhat strong conditions for $f$,
in particular, $f$ has to be smooth.

For the Fermi gas at zero temperature, \cite{2015LS} proved the global well-posedness and orbital stability when $w \in L^1_x \cap L^\I_x$.
\cite{2017CHP} proved global well-posedness when $w$ is the delta measure.
However, we have no results for the {\it asymptotic} stability of the Fermi gas at zero temperature, neither in terms of random field nor in the terms of linear operators. 
Our results in this paper include the asymptotic stability of the Fermi gas at zero temperature.
Moreover, when $d=3$, we have it only in {\it focusing} case, that is, $\hat{w}(0) \le 0$;
this contrasts with the results of \cite{2015LS} and \cite{2017CHP} because their results are in {\it defocusing} case,
that is, $\hat{w}(0) \ge 0$.

\subsection{Main results}
In this paper, we prove the asymptotic stability of $Y=Y_f$ under weak assumptions for $f$ that allows Fermi gas at zero temperature.
Let $Z(t,x,\om)$ be a perturbation from $Y$, that is, $Z(t):= X(t)-Y$.
Then, we have
\begin{displaymath}\label{IVP} \tag{IVP}
	\left\{
	\begin{array}{l} 
		i\pl_t Z = (m - \De + V)Z + VY, \\
		V = w \ast \rh_Z, \\
		\rh_Z = \BE[|Z|^2] + 2\Re \BE[\overline{Y}Z], \\
		Z(0)=Z_0.
	\end{array}
	\right.
\end{displaymath}

Before stating our main results, we define some notations.
For $p \in [1,\I]$ and $\si \in \BR$, we define honogeneous and inhomogeneous Sobolev spaces by
\begin{align*}
	\|u\|_{\dH^{\si,p}} := \||\na|^\si u\|_{L^p},
	\quad \|u\|_{H^{s,p}}:= \|\sd^\si u\|_{L^p}.
\end{align*}
For any $f \in L^2_\xi$, we define $H_f(x):= \CF^{-1}[|f|^2](x)$.
For any radial $f \in L^2_\xi$, we define an even function $h_f \in L^\I(\BR)$ by $h_f(|x|)= H_f(x)$.
For any $f$ and $w$, we define $S(t):= e^{-it(m-\De)}$, where $m:=\wh{w}(0)\|f\|_{L^2_\xi}^2$.
For any $f$ and $w$, we define a linear operator $L=L(f,w)$ by
\begin{align}
	u(t,x) \mapsto L[u](t,x) = 2\Re \BE [\overline{Y(t)}D((w\ast u)Y)(t)], \label{L}
\end{align}
where $D(\ps)(t):= -i \int_0^t S(t-\ta)\ps(\ta)d\ta$.

Our main results are the following.
In these theorems, we assume $\|(1-L(f,w))^{-1}\|_{\CB(L^2_t L^2_x)} < \I$,
but later, we give some sufficient conditions for this to hold.

The first theorem requires the weakest condition for $f$ in this paper.
However, we need somewhat strong conditions for potential $w$,
in particular, we need $\hat{w}(0)=0$ when $d=3$.
We remove it in Theorem \ref{main 3D} by requiring more strong condition for $f$
and more regularity for $Z_0$.
\begin{th.} \label{main}
	Let $d \ge 3$.
	Let $s = \frac{d-3}{2}$.
	Let $f\in L^2_\xi$ be a radial and real-valued function such that $h_f \in L^1(\BR)$.
	Let $w \in L^1(\BR^d)$ be a real-valued function and satisfy
	\begin{align*}
		w \in \dH^{s-\tw,1}_x \cap \dH^{s+\ft,1}_x \cap \dH^{s}_x \cap \dH^{s+\ft}_x.
	\end{align*}
	Assume that $(1-L(f,w))^{-1} \in \CB(\Ltx)$. 
	Then, there exists $\ep_0 > 0$ such that the following holds:
	If $\|Z_0\|_{L^2_\om (\dot{H}^{-s-1/2}_x \cap \dot{H}^{1/4}_x )} \le \ep_0$, 
	then there exists a unique solution to \eqref{IVP}, $Z(t) \in C(\BR, L^2_\om H^\ft_x)$,
	such that $\rh_Z \in L^2_t(\BR, \dH^{-s}_x)$.
	Moreover, $Z(t)$ scatters, that is, 
	there exist $Z_\pm \in L^2_\om H^\ft_x$ and 
	\begin{align*}
		S(-t)Z(t) \to Z_\pm \mbox{ in } L^2_\om H^\ft_x \mbox{ as } t \to \pm \I.
	\end{align*}
\end{th.}

In Theorem \ref{main}, we need $\wh{w}(0) = 0$ when $d=3$.
We can avoid this assumption by requiring stronger condition to $f$ and initial data $Z_0$.
Moreover, we can allow more singular potentials, for example, delta measure.
Define $A_\th[g]$ by
\begin{align} \label{C1}
	A_\th[g]:= \Dk{\int_\BR dv \K{  \int_\BR du \Abs{g(\sqrt{u^2+v^2})|u|^\th} }^2 }^\tw 
\end{align}

\begin{th.} \label{main 3D}
	Let $d=3$.
	Let $w$ be a finite signed Borel measure on $\BR^3$.
	Let $f \in L^2(\BR^3_\xi)$ be a radial and real-valued function.
	Let $f$ satisfy
	$h_f \in L^1(\BR)$ and $A_\th[h_f] < \I$ for some $0 < \th < \tw$.
	Assume that $(1-L(f,w))^{-1} \in \CB(\Ltx)$.
	Then, there exists $\ep_0 > 0$ such that the following holds:
	If $\|Z_0\|_{L^2_\om(\dH^{-1/2}_x \cap \dH^{1/2}_x)} \le \ep_0$,
	then there exists a unique global solution to (IVP), $Z(t) \in C(\BR, L^2_\om \dH^\tw_x)$,
	such that $\rh_Z \in L^2_t(\BR, H^\tw_x)$.
	Moreover, $Z(t)$ scatters, that is, there exist $Z_\pm \in L^2_\om \dH^\tw_x$ such that
	\begin{align*}
		S(-t)Z(t) \to Z_\pm \mbox{ in } L^2_\om \dH^\tw_x \mbox{ as } t \to \pm \I.
	\end{align*}
\end{th.}

\begin{rmk.}
	Note that
	\begin{align*}
		\|Z_0\|_{L^2_\om \dot{H}^{-s-1/2}_x} \lesssim \|Z_0\|_{L^2_\om L^{d'}_x}.
	\end{align*}
	Therefore, for initial data, it is sufficient if $\|Z_0\|_{L^2_\om(L^{d'}_x \cap \dot{H}^{1/4}_x)}$ is small.
	This means, we need only $\ft$-differentiability for initial data in any dimension $d \ge 3$.
	This fact contrasts with the previous studies \cite{2020CS, 2022CS}, which required $(\frac{d}{2}-1)$-differentiability for initial data.
\end{rmk.}

Now, we give some sufficient conditions for $\|(1-L(f,w))^{-1}\|_{\CB(L^2_t L^2_x)}<\I$ to hold.
\begin{prop.} \label{inv 1}
	Let $d \ge 1$.
	Let $w$ be a finite signed Borel measure.
	Let $f \in L^2_\xi$ be real-valued and radial.
	Assume that \eqref{SC} or \eqref{CS} holds:
	\begin{align}
		&\frac{\|\widehat{w}\|_{L^\I_\xi}}{2|\BS^{d-1}|} \int_{\BR^d} \frac{|H_f(x)|}{|x|^{d-2}} dx < 1 \label{SC}, \\
		&\No{\frac{\wh{w}(\xi)}{|\xi|}}_{L^\I_\xi} \int_0^\I |h_f(r)|dr < 1 \label{CS}.
	\end{align}
	Then $(1-L(f,w))^{-1} \in \CB(L^2_t L^2_x)$.
\end{prop.}

Moreover, when $d=3$ and $|f(\xi)|^2=\ch_{\{|\xi| \le 1\}}$, we have the following sufficient condition.
\begin{prop.} \label{cor}
	Let $d=3$ and $|f(\xi)|^2=\ch_{\{|\xi| \le 1\}}$.
	Let $w$ be a finite signed Borel measure on $\BR^3$.
	Let $w$ satisfy $\hat{w}$ is real-valued and
	\ALN{
		-\de \le \hat{w}(\xi) \le \frac{\de_0}{\Br{\log|\xi|}}
	}
	for sufficiently small absolute constants $\de, \de_0 > 0$.
	Then, we have $(1-L(f,w))^{-1} \in \CB(L^2_t L^2_x)$.
\end{prop.}

\begin{rmk.}
	Our all results include the case $|f(\xi)|^2=\ch_{\{|\xi| \le 1\}}$, that is, Fermi gas at zero temperature.
	By Lemma \ref{decay rate}, we have $|h_f(r)| \lesssim \Br{r}^{-\frac{d+1}{2}}$.
	Thus, we have $h_f \in L^1(\BR)$ when $d \ge 2$.
	Moreover, the decay rate $|H_f(x)| \le \bx^{-\frac{d+1}{2}}$ implies
	\begin{align*}
		&\frac{\|\widehat{w}\|_{L^\I_\xi}}{2|\BS^{d-1}|} \int_{\BR^d} \frac{|H_f(x)|}{|x|^{d-2}} dx < \I
	\end{align*}
holds when $d \ge 4$, and 
\begin{align*}
\int_0^\I |h_f(r)|dr < \I
\end{align*}
holds when $d \ge 2$.

Let $d=3$.
Then we have
	\begin{align*}
		|h_f(\sqrt{u^2+v^2})|
		&\lesssim \frac{1}{1+u^2+v^2} \\
		&\le \frac{1}{(1+u^2)^{11/16}(1+v^2)^{5/16}}.
	\end{align*}
	Therefore, we have for $\th = \frac{1}{4}$
	\begin{align*}
		A_\th[h_f]^2
		\lesssim \int_\BR \frac{1}{(1+v^2)^{5/8}} dv \K{ \int_\BR \frac{|u|^{1/4}}{(1+u^2)^{11/16}} du }^2
		< \I.
	\end{align*}
\end{rmk.}

\begin{rmk.}
	When $d=3$ and $|f(\xi)|^2 = \ch_{\{|\xi| \le 1\}}$,
	we have the scattering results only in focusing case.
	This is surprising, but we have no better explanation for this reason than that $m_f(\tau,\xi)$ is upper bounded when $d = 3$.
	(See \eqref{m_f formula}.)
\end{rmk.}

\subsection{Strategy}
Roughly speaking, we reduce \eqref{IVP}, a system for two unknown function $Z$ and $\rh=\rh_Z$, to an equation for one unknown function $\rh$.
This idea is inspired by \cite{2015LS, 2014LS, 2018CHP}.
We can see it as a random field version of the density functional theory.
However, it contrasts with \cite{2020CS, 2022CS} since they considered the system of $Z$ and $\rh$.
 
Before starting to explain our strategy, we define some notations.
Let $S_V(t,t_0)$ be the propagator of the following linear Schr\"{o}dinger equation with time-dependent potential:
\begin{align} \label{VLS}
	(i\pl_t - m + \De-V(t,x))u =0,
	\quad u(t,x):\BR \times \BR^d \to \BC,
\end{align}
where $V(t,x):\BR\times \BR^d \to \BR$.
Namely, $u(t) := S_V(t,t_0)u(t_0)$ is the solution to \eqref{VLS}.
Note that $S(t)= S_0(t)$.
Define $D_V$ by
\begin{align*}
	D_V(\ps)(t):= -i\int_0^t S_V(t,\ta) \ps(\ta) d\ta.
\end{align*}
Note that $D = D_0$.
\eqref{IVP} is equivalent to the following:
\begin{displaymath}
	\left\{
	\begin{array}{l} 
		Z(t) = S_V(t)Z_0 + D_V(VY)(t), \\
		V = w \ast \rh, \\
		\rh = \BE[|Z|^2] + 2\Re \BE[\overline{Y}Z].
	\end{array}
	\right.
\end{displaymath}
Therefore, we get the equation for $\rh$:
\begin{align*}  \tag{IVP*} \label{IVP*} 
    \rh = \BA(\rh),
\end{align*}
where
\begin{align*}
    &\BA(\rh) := \BE[|S_V(t)Z_0|^2] + 2\Re \BE[(S_V(t)Z_0)\overline{D_V(VY)}] +\BE[|D_V(VY)|^2]\\
	&\qquad \qquad + 2\Re \BE[\overline{Y}S_V(t)Z_0] + 2\Re \BE [\overline{Y}D_V(VY)], \\
	&V = w \ast \rh.
\end{align*}

The goal of this paper is to prove the following theorems:
\begin{th.}\label{rh}
	Under the same assumptions as in Theorem \ref{main},
	there exists $\ep_0 > 0$ such that the following holds: 
	For any $\|Z_0\|_{L^2_\om(\dH^{-s-1/2}_x \cap \dH^{1/4}_x)} \le \ep_0$,
	there exists a unique solution $\rh \in L^2_t(\BR,\dH^{-s}_x)$ to \eqref{IVP*}.
\end{th.}

\begin{th.}\label{rh 3D}
	Under the same assumptions as in Theorem \ref{main 3D},
	there exists $\ep_0 > 0$ such that the following holds: 
	For any $\|Z_0\|_{L^2_\om(\dH^{-1/2}_x \cap \dH^{1/2}_x)} \le \ep_0$,
	there exists a unique solution $\rh \in L^2_t(\BR,H^\tw_x)$ to \eqref{IVP*}.
\end{th.}

The proof of the main results are reduced to proving the above theorems
because Theorem \ref{rh}, \ref{rh 3D} and the following lemmas imply the main results.
We prove these main lemmas in Sect 7.
\begin{lem.} \label{scattaring lemma}
	Let $d \ge 3$ and $0 \le \si \le \frac{d}{2}-1$.
	Let $Z_0 \in L^2_\om H^\si_x$ and $h_f \in L^1(\BR)$.
	Assume that
	\begin{align*}
		V \in L^2_t(\BR,H^{\si}_x \cap H^{\si,\I}_x \cap \dH^{-\tw}_x).
	\end{align*}
	Define $Z$ by $Z(t) := S_V(t)Z_0 + D_V(VY)(t)$.
	Then $Z(t) \in C(\BR, L^2_\om H^{\si}_x)$ and $Z(t)$ scatters in the following sense:
	There exist $Z_\pm \in L^2_\om H^{\si}_x$ such that
	\begin{align*}
		S(-t)Z(t) \to Z_\pm \;\; \mbox{in} \;\; L^2_\om H^\si_x \;\; \mbox{as} \;\; t \to \pm \I.
	\end{align*}
\end{lem.}

\begin{lem.} \label{scattaring lemma 1.5}
	Let $d = 3$.
	Let $Z_0 \in L^2_\om H^\tw_x$ and $h_f \in L^1(\BR)$.
	Assume that
	\begin{align*}
		V \in L^2_t(\BR,H^\tw_x).
	\end{align*}
	Define $Z$ by $Z(t) := S_V(t)Z_0 + D_V(VY)(t)$.
	Then $Z(t) \in C(\BR, L^2_\om \dH^\tw_x)$ and $Z(t)$ scatters in the following sense:
	There exist $Z_\pm \in L^2_\om \dH^\tw_x$ such that
	\begin{align*}
		S(-t)Z(t) \to Z_\pm \mbox{ in } L^2_\om \dH^\tw_x \mbox{ as } t \to \pm \I.
	\end{align*}
\end{lem.}

To prove Theorems \ref{rh} and \ref{rh 3D}, we transform \eqref{IVP*}.
We expand the last term to the following form:
\begin{align*}
	2\Re \BE [\overline{Y}D_V(VY)] = L(\rh)  + Q(w \ast \rh),
\end{align*}
where
\begin{align}
	&L(\rh) = 2\Re \BE [\overline{Y}D((w\ast \rh)Y)]\\
	&Q(V) = -2 \Re \BE \Dk{ \ov{Y(t)} \int_0^t S(t-s) V(s) \int_0^s S_V(s,\ta) (V(\ta)Y(\ta))d\ta ds} \label{Q}.
\end{align}
If $(1-L)^{-1} \in \CB(L^2_{t,x})$, then \eqref{IVP*} is equivalent to the following equation:
\begin{align} \label{rh2}
	\rh = (1-L)^{-1} ( \BA_1(\rh) + \cdots + \BA_5(\rh) ),
\end{align}
where
\begin{align*}
	&\BA_1(\rh) = \BE[|S_V(t)Z_0|^2], \\
	&\BA_2(\rh) = 2\Re \BE[(S_V(t)Z_0)\overline{D_V(VY)}], \\
	&\BA_3(\rh) = \BE[|D_V(VY)|^2], \\
	&\BA_4(\rh) = 2\Re \BE[\overline{Y}S_V(t)Z_0],  \\
	&\BA_5 (\rh) = Q(V) =-2\Re \BE \Dk{ \ov{Y(t)} \int_0^t S(t-s) V(s) \int_0^s S_V(s,\ta) (V(\ta)Y(\ta))d\ta ds}.
\end{align*}
The rest of this paper will be devoted to finding a fixed point of $\Ph[\rh]:=(1-L)^{-1} ( \BA_1(\rh) + \cdots + \BA_5(\rh) )$.
To prove Theorem \ref{rh 3D}, we need cancellation of quadratic terms used in \cite{2020CS}, but we explain it in Sect 5.2.

\subsection{Organization and acknowledgements}
This paper is organized as follows.
In Sect 2, we give some sufficient conditions for $1-L$ to be invertible.
In Sect 3, we prove the Strichartz estimates for the propagator $S_V(t)$.
In Sect 4, we generalize the randomized Strichartz estimates proved in \cite{2022CS} to the propagator $S_V(t)$.
These estimates are the key estimates in this paper.
In Sect 5, we collect the necessary estimates to prove main result as applications of the key estimates.
Moreover, in Sect 5.2, we explain the strategy of the proof of Theorem \ref{rh 3D}.
In Sect 6, we prove the main theorems.
In Sect 7, we prove main lemmas, that is, Lemmas \ref{scattaring lemma} and \ref{scattaring lemma 1.5}.

\subsection*{\centerline{Acknowledgements}}
The author would like to thank Prof. Kenji Nakanishi in Kyoto University.
He gave the author many useful comments throughout this paper.
The author also would like to thank Prof. Mikihiro Fujii in Kyushu University for his contribution to Proposition \ref{upper bound}.

\section{Linear response}
In this section, we discuss the invertibility of $1-L$.
It is known that $L$ is a space-time Fourier multiplier with symbol $\widehat{w}(\xi)m_f(\ta,\xi)$, where
\begin{align} \label{m_f}
	m_f(\ta,\xi) = -2\int_0^\I e^{-i\ta t} \sin(|\xi|^2 t) H_f(2\xi t) dt.
\end{align}
For this fact, see by \cite[Lemma 5.6]{2020CS}; and see also \cite[Proposition 1]{2014LS}.
Therefore, $(1-L)^{-1} \in \CB(L^2_t H^{-s}_x)$ is equivalent to $(1-L)^{-1} \in \CB(L^2_t L^2_x)$.
Some sufficient conditions for $(1-L)^{-1} \in \CB(\Ltx)$ are known.
For example, see \cite[Corollary 1]{2014LS}.
The author could not find any clear statement in the literature, but \eqref{m_f} immediately implies Proposition \ref{inv 1}

In the last of this section, we prove the invertibility of $1-L$ when $|f(\xi)|^2=\ch_{\{|\xi| \le 1\}}$, that is, the case of Fermi gas at zero temperature.

First, we have the following decay rates of $H_f$.
\begin{lem.} \label{decay rate}
	Let $d\ge 1$ and $|f(\xi)|^2=\ch_{\{|\xi|\le 1\}}$. 
	Then we have
	\begin{align}
		|H_f(x)| \lesssim \bx^{-\frac{d+1}{2}}.
	\end{align}
\end{lem.}

\begin{proof}
	It is easy to verify that
	\ALN{
		H_f(x) &= C_d \int_{-1}^1 (1-t^2)^{\frac{d-1}{2}} \cos(|x|t) dt.
	}
	Here, we have Poisson's integral representations for the Bessel function $J_\nu$, that is, 
	\ALN{
		J_\nu(x) = C_\nu x^\nu \int_{-1}^1 (1-t^2)^{\nu-1/2} \cos(xt) dt,  
	}
	for $\nu > -1/2$ and $x>0$.
	See \cite[Appendix B]{2014G} for the Bessel functions.
	Hence, we have
	\ALN{
		H_f(x) = C_d |x|^{-d/2} J_{d/2}(|x|).
	}
	Here, we have the decay rate of the Bessel functions:
	\ALN{
		|J_\nu(x)| \leq  \frac{C_\nu}{|x|^{1/2}}.
	}
	Since $H_f$ is bounded, it follows from the above argument that
	\AL{
		|H_f(x)| \leq C_d \Br{x}^{-\frac{d+1}{2}}.
	}
\end{proof}

Moreover, when $d=3$, we have the following formula of $m_f$.
\begin{prop.} \label{upper bound}
	Let $d=3$ and $|f(\xi)|^2 = \ch_{\{|\xi| \le 1\}}$.
	For $(\ta,\xi) \in \BR^{1+3}$, define $\ep = \ep(\ta,\xi)$ by
	\ALN{
		\ep := \min \K{ \Abs{2-\frac{\ta}{|\xi|}}, \Abs{2+\frac{\ta}{|\xi|}}}.
	}
	Then, there exists $b \in L^\I(\BR^{1+3})$ such that
	\AL{ \label{m_f formula}
		m_f(\ta,\xi) = \frac{1}{2\sqrt{2\pi}} \min(\log[\max(\ep, |\xi|)], 0) + b(\ta,\xi).
	} 
	In particular, $\Re[m_f(\tau,\xi)]$ is upper bounded on $\BR^{1+3}$.
\end{prop.}

The formula \eqref{m_f formula} immediately implies Proposition \ref{cor}

\begin{proof}
	It follows from the direct calculation that
	\ALN{
		H_f(x) = \sqrt{\frac{2}{\pi}} \frac{1}{|x|^2} \K{\frac{\sin(|x|)}{|x|} - \cos(|x|)}.
	}
	
	Let $\Om := [-1,1]\times [-1,1]^3$.
	It is easy to see that $m_f(\tau, \xi)$ is bounded on $\BR^{1+3}\setminus \Om$.
	Thus, we consider the case $(\ta,\xi) \in \Om$.
	Let $k:= |\xi|$.
	Then we have
	\ALN{
		I(\tau,k) &:= m_f(\tau, \xi) \\
		&= \frac{2\sqrt{2}}{\sqrt{\pi}} \int_0^\I e^{-i\tau t} \sin(tk^2) \frac{1}{(2tk)^2}\K{\frac{\sin(2tk)}{2tk}-\cos(2tk)} dt \\
		&= \frac{1}{\sqrt{2\pi}} \int_0^\I e^{-i\frac{\tau}{k} t} \sin(tk) \frac{1}{t^2}\K{\frac{\sin(2t)}{2t}-\cos(2t)} \frac{dt}{k} \\
		&= \frac{1}{\sqrt{2\pi}} \int_0^1 e^{-i\frac{\tau}{k} t} \sin(tk) \frac{1}{t^2}\K{\frac{\sin(2t)}{2t}-\cos(2t)} \frac{dt}{k} \\
		& \quad + \frac{1}{\sqrt{2\pi}} \int_1^\I e^{-i\frac{\tau}{k} t} \sin(tk) \frac{\sin(2t)}{2t^3} \frac{dt}{k}
		- \frac{1}{\sqrt{2\pi}} \int_1^\I e^{-i\frac{\tau}{k} t} \sin(tk) \frac{\cos(2t)}{t^2}  \frac{dt}{k} \\
		&=:I_1(\tau,k) + I_2(\tau,k) + I_3(\tau,k).
	}
	We have
	\ALN{
		|I_1(\tau,k)| &\lesssim \int_0^1 t dt \lesssim 1, \\
		|I_2(\tau,k)| &\lesssim \int_1^\I \frac{1}{t^2} dt \lesssim 1.
	}
	Next we estimate $I_3(\tau,k)$.
	Integrating by parts with respect to $\cos(2t)$ twice, we have
	\ALN{
		I_3 \leq \frac{1+\frac{|\tau|}{k}}{\Abs{1-\frac{\tau^2}{4k^2}}}.
	}
	This implies $m_f(\tau,\xi)$ is bounded on $\Om \cap \{\ep >1 \}$.
	Let $(\ta,\xi) \in \Om$ and $\ep \le 1$.
	Then, we can write $\tau/k = \pm(2+\ep_0)$ with $\ep_0 \in [-1,1]$.
	Note that $|\ep_0| = \ep$.
	We have
	\ALN{
		\Re[I_3(\tau,k)] &= - \frac{1}{\sqrt{2\pi}} \int_1^\I \cos\K{\frac{\tau}{k}t}\sin(tk) \frac{\cos(2t)}{t^2} \frac{dt}{k} \\
		&= - \frac{1}{\sqrt{2\pi}} \int_1^\I \cos((2+\ep_0)t)\sin(tk) \frac{\cos(2t)}{t^2} \frac{dt}{k} \\
		&= -\frac{1}{2\sqrt{2\pi}} \int_1^\I \cos(\ep_0 t)\sin(tk) \frac{1}{t^2} \frac{dt}{k} 
		-\frac{1}{2\sqrt{2\pi}} \int_1^\I \cos((4 +\ep_0) t)\sin(tk) \frac{1}{t^2} \frac{dt}{k} \\
		&=: \frac{1}{2\sqrt{2\pi}}(J_1(\tau,k) + J_2(\tau,k)).
	}
	An integration by parts with respect to $\cos((4+\ep_0)t)$ and $\ep_0 \in [-1,1]$ imply $J_2(\tau,k)$ is bounded.
	Next, we estimate $J_1(\tau,k)$.
	We have
	\ALN{
		J_1(\tau,k) &= - \int_1^{1/k} \cos(\ep t)\sin(tk) \frac{1}{t^2} \frac{dt}{k} - \int_{1/k}^\I \cos(\ep t)\sin(tk) \frac{1}{t^2} \frac{dt}{k}\\
		&=: J_1^1(\tau,k) + J_1^2(\tau,k).
	}
	It follows from changing variables that $J_1^2(\tau,k)$ is bounded.
	First, we assume that $k \ge \ep$. 
	Then, we have
	\AL{
		J_1^1(\ta,k) &= - \int_1^{1/k} \frac{dt}{t} + \int_1^{1/k} \K{1-\cos(\ep t)} \frac{dt}{t} + \int_1^{1/k} \cos(\ep t) \K{\frac{1}{t}-\frac{\sin(kt)}{kt^2} } dt \label{zurasu} \\
		&=\log|k| + (\mbox{bounded term}). \nonumber
	}
	Next we assume that $k < \ep$.
	Then, we have
	\ALN{
		J_1^1(\ta,k) &= - \int_0^{1/\ep} \cos(\ep t)\sin(tk) \frac{1}{t^2} \frac{dt}{k} - \int_{1/\ep}^{1/k} \cos(\ep t)\sin(tk) \frac{1}{t^2} \frac{dt}{k}.
	}
	By integration by parts, we have the second term is bounded.
	For the first term, by the same way as \eqref{zurasu}, we have
	\ALN{
		- \int_0^{1/\ep} \cos(\ep t)\sin(tk) \frac{1}{t^2} \frac{dt}{k} = \log \ep + (\mbox{bounded term}).
	}
	From the above, we have
	\ALN{
		J_1^1(\ta,k) = \log[\max(\ep,k)] + (\mbox{bounded term}).
	}
	
	Next, we show that $\mathrm{Im} [I_3(\ta,k)]$ is bounded.
	In the following calculation, we ignore any constants, that is, we write $A \sim B$ if there exists $C \in \BC$ such that $A = CB$.
	We can assume that $\ta/k$ is near $2$.
	\ALN{
		\mathrm{Im}[I_3(\ta,k)] &\sim \int_1^\I \sin \K{\frac{\tau}{k}t}\sin(tk) \frac{\cos(2t)}{t^2} \frac{dt}{k} \\
		&\sim \int_1^\I \sin \K{(2+\ep_0)t}\sin(tk) \frac{\cos(2t)}{t^2} \frac{dt}{k} \\
		&\sim \int_1^\I \sin(\ep_0 t) \sin(kt) \frac{dt}{kt^2} + \int_1^\I \sin((4+\ep_0)t)\sin(kt) \frac{dt}{kt^2} \\
		&=: J_3(\ta,k) + J_4(\ta,k).
	}
	By integration by parts, we have $J_4(\ta,k)$ is bounded.
	We have
	\ALN{
		J_3(\ta,k) &= \int_1^{1/k} \sin(\ep_0 t) \sin(kt) \frac{dt}{kt^2} + \int_{1/k}^\I \sin(\ep_0 t) \sin(kt) \frac{dt}{kt^2} \\
		&=: J_3^1(\ta,k) + J_3^2(\ta,k).
	}
	By changing variables, we have $J_3^2(\ta,k)$ is bounded.
	When $\ep = |\ep_0| \le k$, we have
	\ALN{
		|J_3^1(\ta,k)| &\le \int_1^{1/k} (\ep t)(kt) \frac{dt}{kt^2} \\
		&\le \frac{\ep}{k} \\
		&\le 1.
	}
	Let $\ep > k$.
	We have
	\ALN{
		J_3^1(\ta,k) = \int_1^{1/\ep}\sin(\ep_0 t) \sin(kt) \frac{dt}{kt^2} + \int_{1/\ep}^{1/k}\sin(\ep_0 t) \sin(kt) \frac{dt}{kt^2}
	}
	By changing variables, we have the first term is bounded.
	By integration by parts, we have the second term is bounded.
\end{proof}

\section{Preliminaries}
We prove the following Strichartz estimates for $S_V(t,\ta)$.
\begin{lem.} \label{SV}
	Let $d \ge 2$.
	Let $p_0, q_0,p_1,q_1 \in [2,\I]$ satisfy $\frac{2}{p_n}+\frac{d}{q_n}=\frac{d}{2}$ and $(p_n,q_n,d)\neq (2,\I,2)$ for $n=0,1$
	and $(p_0,q_0,d) \neq (\I,2,2)$.
	Let $\frac{1}{c} = \frac{\si+1}{d}$.
	Let $\si$ satisfy
	\begin{align*}
		0\le \si \le \frac{d}{2}-1, \quad \frac{\si}{d}< \frac{1}{q_0} .
	\end{align*}
	Then there exists a monotone increasing function $\ph:[0,\I)\to[0,\I)$ such that the following hold:
	\begin{align}
		&\||\na|^\si S_V(t)u_0\|_{L^{p_0}_t L^{q_0}_x}
		\le \ph(\|V\|_{L^2_t \dot{H}^{\si, c}_x})
		\||\na|^\si u_0\|_{L^2_x}, \label{homo} \\
		&\No{|\na|^\si \int_\BR S_V(t,\ta) u(\ta) d\ta}_{L^{p_0}_t L^{q_0}_x}
		\le \ph(\|V\|_{L^2_t \dH^{\si,c}_x})
		\||\na|^\si u\|_{L^{p_1'}_t L^{q_1'}_x}, \label{inhomo}\\
		&\No{|\na|^\si \int_0^t S_V(t,\ta) u(\ta) d\ta}_{L^{p_0}_t L^{q_0}_x}
		\le \ph(\|V\|_{L^2_t \dH^{\si,c}_x})
		\||\na|^\si u\|_{L^{p_1'}_t L^{q_1'}_x}. \label{retarted}
	\end{align}
\end{lem.}

\begin{rmk.}
	In this paper, we denote any monotone increasing function from $[0, \I)$ to $[0,\I)$ by $\ph$.
	This is like the constant $C$.
	We do not distinguish a monotone increasing function from another one.
\end{rmk.}

\begin{proof}
	Let $u_0 \in \CS(\BR^d)$.
	First, we have
	\begin{align*}
		|\na|^\si S_V(t) u_0 &= S(t)|\na|^\si u_0 - i \int_0^t S(t-\ta) |\na|^\si (V(\ta) S_V(\ta)u_0 ) d\ta.
	\end{align*}
	We define $\th_0, \mu_0, \mu_1 \in [2,\I]$ by
	\begin{align*}
		\frac{1}{\th_0'}= \frac{1}{p_0}+ \frac{1}{2},
		\quad \frac{2}{\th_0} + \frac{d}{\mu_0} = \frac{d}{2},
		\quad  \frac{1}{\mu_1} = \frac{1}{q_0} - \frac{\si}{d}.
	\end{align*}
	Note that we have
	\begin{align}
		\quad \frac{1}{\mu_0'} = \frac{1}{c}+\frac{1}{\mu_1} = \frac{1}{q_0}+\frac{1}{d}. \label{expo}
	\end{align}
	Then we have by the fractional Leibniz rule
	\begin{align*}
		&\||\na|^\si S_V(t) u_0\|_{L^{p_0}_t L^{q_0}_x}
		\le
		C_0 \||\na|^\si u_0\|_{L^2_x}
		+ C_0 \||\na|^\si(V(t)S_V(t)u_0)\|_{L^{\th_0'}_t L^{\mu_0'}_x} \\
		&\quad \le C_0 \||\na|^\si u_0\|_{L^2_x} + C_1 \||\na|^\si V\|_{L^2_t L^c_x}\|S_V(t)u_0\|_{L_t^{p_0} L_x^{\mu_1}} 
		+ C_1 \|V\|_{L^2_t L^d_x} \||\na|^\si S_V(t)u_0\|_{L^{p_0}_t L^{q_0}_x} \\
		&\quad \le C_0 \||\na|^\si u_0\|_{L^2_x}
		+ C_1\|V\|_{L^2_t \dH^{\si,c}_x}
		  \||\na|^\si S_V(t)u_0\|_{L^{p_0}_t L^{q_0}_x}.
	\end{align*}
	Note that $C_0 = C(d,p_0)$ and $C_1=C(d,p_0,\si)$.
	Therefore, if $C_1\|V\|_{L^2_t \dH^{\si,c}_x} \le \tw$, we have
	\begin{align*}
		\||\na|^\si S_V(t)u_0\|_{L^{p_0}_t L^{q_0}_x} \le 2C_0 \||\na|^\si u_0\|_{L^2_x}.
	\end{align*}
	Now we assume that $\|V\|_{L^{p_0}_t(\BR,\dH^{\si,c}_x)} \ge \frac{1}{2C_1} = :\de$.
	There exist intervals $I_n \subset \BR$ for $n=1, \dots , N+1$ such that $\BR=\cup_{n=1}^{N+1} I_n$ and 
	\begin{align}
		&\frac{\de}{2} \le \|V\|_{L^{p_0}_t(I_n, \dH^{\si,c}_x)} \le \de, \mbox{ for } 1 \le n \le N; \label{N}\\
		&\|V\|_{L^{p_0}_t(I_{N+1}, \dH^{\si,c}_x)} \le \de. \nonumber
	\end{align}
	By the same argument as above, we have
	\begin{align*}
		\||\na|^\si S_V(t)u_0\|_{L^{p_0}_t(I_n, L^{q_0}_x)} \le 2C_0 \||\na|^\si u_0\|_{L^2_x}.
	\end{align*}
	Therefore, we have
	\begin{align*}
		\||\na|^\si S_V(t)u_0\|_{L^{p_0}_t(\BR, L^{q_0}_x)} \le 2(N+1)C_0 \||\na|^\si u_0\|_{L^2_x}
	\end{align*}
	We have by \eqref{N} that
	\begin{align*}
		N \le \K{2/\de}^{p_0} \|V\|_{L^{p_0}_t(\BR, \dH^{\si,c}_x)}^{p_0}.
	\end{align*}
	Therefore, we have
	\begin{align*}
		\||\na|^\si S_V(t)u_0\|_{L^{p_0}_t(\BR, L^{q_0}_x)}
		\le 2\K{(2/\de)^{p_0} \|V\|_{L^{p_0}_t(\BR, \dH^{\si,c}_x)}^{p_0}+1} C_0 \||\na|^\si u_0\|_{L^2_x}.
	\end{align*}

	Next, we consider \eqref{inhomo} and \eqref{retarted}.
	We only prove \eqref{retarted}; since the other can be shown in the same way.
	Note that 
	\begin{align*}
		S_V(t,\ta) = S(t-\ta) - i\int_\ta^t S(t-s)V(s)S_V(s,\ta)ds.
	\end{align*}
	Hence, we have
	\begin{align*}
		&\No{|\na|^\si \int_0^t S_V(t,\ta)u(\ta)d\ta}_{L^{p_0}_t L^{q_0}_x} \\
		&\quad \le \No{|\na|^\si \int_0^t S(t-\ta)u(\ta)d\ta}_{L^{p_0}_t L^{q_0}_x}
		+ \No{\int_0^t \int_\ta^t S(t-s)|\na|^\si V(s) S_V(s,\ta)u(\ta)dsd\ta}_{L^{p_0}_t L^{q_0}_x} \\
		&\quad \le C_0 \||\na|^\si u\|_{L^{p_1'}_t L^{q_1'}_x} + \No{\int_0^t S(t-s)|\na|^\si V(s) \int_0^s S_V(s,\ta)u(\ta)d\ta ds}_{L^{p_0}_t L^{q_0}_x}.
	\end{align*}
	For the same exponents as \eqref{expo}, we have
	\begin{align*}
		&\No{|\na|^\si \int_0^t S_V(t,\ta)u(\ta)d\ta}_{L^{p_0}_t L^{q_0}_x} \\
		&\quad \le C_0 \||\na|^\si u\|_{L^{p_1'}_t L^{q_1'}_x} + C_0 \No{|\na|^\si V(t) \int_0^t S_V(t,\ta)u(\ta)d\ta}_{L^{\th_0'}_t L^{\mu_0'}_x} \\
		&\quad \le C_0 \||\na|^\si u\|_{L^{p_1'}_t L^{q_1'}_x}
		+ C_1 \||\na|^\si V\|_{L^2_t L^c_x}	\No{\int_0^t S_V(t,\ta)u(\ta)d\ta}_{L^{p_0}_t L^{\mu_1}_x} \\
		&\qquad + C_1 \|V\|_{L^2_t L^d_x}
		\No{|\na|^\si \int_0^t S_V(t,\ta)u(\ta)d\ta}_{L^{p_0}_t L^{q_0}_x} \\
		&\quad \le C_0 \||\na|^\si u\|_{L^{p_1'}_t L^{q_1'}_x}
		+ C_2 \|V\|_{L^2_t\dH^{\si,c}_x} \No{|\na|^\si \int_0^t S_V(t,\ta)u(\ta)d\ta}_{L^{p_0}_t L^{q_0}_x}.
	\end{align*}
	Therefore, if $C_2 \|V\|_{L^2_t\dH^{\si,r}_x} \le \tw$, we have
	\begin{align}
		\No{|\na|^\si \int_0^t S_V(t,\ta) u(\ta) d\ta}_{L^{p_0}_t L^{q_0}_x}
		\le 2C_0 \||\na|^\si u\|_{L^{p_1'}_t L^{q_1'}_x}. 
	\end{align}
	When $C_2 \|V\|_{L^2_t \dH^{\si,c}_x} \ge \tw$, we can prove \eqref{retarted} by the same argument as \eqref{homo}.
\end{proof}

\section{Randomized Strichartz estimates}
Define $T:Z_0(\om,x) \to T(Z_0)(t,x)$ by
\begin{align*}
	T(Z_0)(t,x) := \BE[\overline{Y(t,x)} (S(t)Z_0)(x)].
\end{align*}
The (formal) adjoint operator of $T$ is
\begin{align*}
	T^{*}(\ps)(\om,x) = \int_{I} S(\ta)^* (Y(\ta)\ps(\ta)) d\ta.\textbf{}
\end{align*}

The following estimates are known:
\begin{prop.}[Proposition 4.7 \cite{2022CS}] \label{RS}
	Let $d \ge 1$ and $\si \in \BR$.
	Then there exist a constant that is independent of $I \subset \BR$ and the following estimates hold:
	\begin{align}
		&\|T(Z_0)\|_{L^2_t(I, \dot{H}^\si_x)} \le C \|h_f\|_{L^1(\BR)} \|Z_0\|_{L^2_\om \dot{H}^{\si-1/2}_x}, \label{T} \\
		&\|T^*(\ps)\|_{L^2_\om \dot{H}^\si_x} \le C \|h_f\|_{L^1(\BR)} \|\ps\|_{L^2_t(I,\dot{H}^{\si-1/2}_x)}. \label{T^*}
	\end{align}
\end{prop.}
By Christ-Kiselev's lemma and Proposition \ref{RS}, we have
\begin{prop.}[Proposition 4.7 \cite{2022CS}] \label{MA}
	Let $d \ge 1$ and $\si \in \BR$.
	Let $p,q \in [2,\I]$ satisfy $\frac{2}{p}+\frac{d}{q} = \frac{d}{2}$ and $p >2$.
	Then we have
	\begin{align*}
		\No{ |\na|^\si \BE\Dk{\ov{Y(t)}\int_0^t S(t-\ta) \Ps(\ta )d\ta} }_\Ltx 
		\lesssim \|h_f\|_{L^1(\BR)} \||\na|^{\si-1/2} \Ps\|_{L^{p'}_t L^2_\om L^{q'}_x }.
	\end{align*}
\end{prop.}

Now we define generalized operators.
\begin{def.}
	Define the operators $T_V$ and $T_V^*$ by
	\begin{align*}
		&T_V(Z_0)(t,x) := \BE[\overline{Y(t,x)}(S_V(t)Z_0)(x)], \\
		&T_V^*(\ps)(\om,x) := \int_\BR S_V(\ta)^*(Y(\ta) \ps(\ta))(x) d\ta.
	\end{align*}
\end{def.}

The following proposition is the main tool of our analysis.
\begin{prop.} \label{fundamental lemma}
	Let $d \ge 2$.
	For $\si \in \BR$, let $\frac{1}{c} := \frac{\si + 1}{d}$ and $\frac{1}{\wt{c}} := \frac{\si + 1/2}{d}$.
	Let $f\in L^2_\xi$ satisfy $h_f \in L^1(\BR)$.
	Then the following estimates hold.
	\begin{itemize}
		\item If $-\frac{d}{2} + \frac{3}{2} \le \si \le \frac{d}{2}-\tw$,
		then there exists a monotone increasing function $\ph:[0,\I) \to [0,\I)$ such that
	           \begin{align} \label{TV+}
	           	\|T_V(Z_0)\|_{L^2_t \dot{H}^{\si}_x}
	           	\le \|h_f\|_{L^1(\BR)}
	           	\ph(\|V\|_{L^2_t \dH^{|\si-1/2|,\wt{c}}_x})
	           	\|Z_0\|_{L^2_\om \dot{H}^{\si-1/2}_x }.
	           \end{align}

		\item If $-\frac{d}{2}+1 \le \si \le \frac{d}{2}-1$, 
		then there exists a monotone increasing function $\ph:[0,\I) \to [0,\I)$ such that
		      \begin{align}\label{TV*+}
		      	\No{ T_V^*(\ps) }_{L^2_\om \dot{H}^{\si}_x}
		      	\le \|h_f\|_{L^1(\BR)}
		      	\ph(\|V\|_{L^2_t \dH^{|\si|,c}_x})
		      	\|\ps\|_{L^2_t \dot{H}^{\si-1/2}_x}.
		      \end{align}
	\end{itemize}
\end{prop.}

\begin{proof}
	First, we consider \eqref{TV*+} when $\si \ge 0$.
	$S_V(t)$ satisfies the following equation:
	\begin{align}
		&S_V(\ta) = S(\ta) - i \int_0^\ta S(\ta-r) V(r) S_V(r) dr, \label{S_V}\\
		&S_V(\ta)^* = S(-\ta) + i \int_0^\ta S_V(r)^* V(r) S(r-\ta) dr. \label{S_V^*}  
	\end{align}
	Therefore, we have
	\begin{align*}
		\|T_V^*(\ps)\|_{L^2_\om \dH^\si_x}
		&=\No{|\na|^\si \int_\BR S_V(\ta)^* Y(\ta) \ps(\ta) d\ta }_{L^2_\om L^2_x} \\ 
		&\le \No{|\na|^\si \int_\BR S(-\ta) Y(\ta) \ps(\ta) d\ta }_{L^2_\om L^2_x} \\
		&\quad + \No{|\na|^\si  \int_\BR \int_0^\ta S_V(r)^* V(r) S(r-\ta) dr Y(\ta) \ps(\ta) d\ta }_{L^2_\om L^2_x} \\
		&=:A + B.
	\end{align*}
	We have by Proposition \ref{RS}
	\begin{align*}
		A \lesssim \|h_f\|_{L^1(\BR)} \||\na|^{\si-\tw} \ps\|_\Ltx.
	\end{align*}
	Next, we estimate $B$.
	We have
	\begin{align*}
		B
		&\le \No{ |\na|^\si  \int_0^\I d\ta \int_0^\ta dr S_V(r)^* V(r) S(r-\ta) Y_f(\ta) \ps(\ta) }_{L^2_\om L^2_x} \\
		&\qquad + \No{|\na|^\si \int_{-\I}^0 d\ta \int_0^\ta dr S_V(r)^* V(r) S(r-\ta) Y_f(\ta) \ps(\ta) }_{L^2_\om L^2_x} \\
		&=: B_1 + B_2.
	\end{align*}
	We only estimate $B_1$ since $B_2$ can be estimated in the same way.
	Let
	\begin{align} \label{qrc}
		\frac{1}{q}=\frac{1}{2}-\frac{1}{2d},
		\quad \frac{1}{q'}=\frac{1}{c}+\frac{1}{\mu}.
	\end{align}
	By Lemma \ref{SV}, we have
	\begin{align*}
		B_1
		&=
		\No{ |\na|^\si \int_0^\I S_V(r)^* V(r) S(r) \int_r^\I S(-\ta) Y_f(\ta) \ps(\ta) d\ta dr }_{L^2_\om L^2_x} \\
		&\le
		\ph(\|V\|_{L^2_t \dH^{\si,c}_x})
		\No{ |\na|^\si V(t) S(t) \int_t^\I S(-\ta) Y(\ta) \ps(\ta) d\ta}_{L^2_\om L^{4/3}_t L^{q'}_x}\\
		&\le
		\ph(\|V\|_{L^2_t \dH^{\si, c}_x})
		\||\na|^\si V\|_{L^2_t L^c_x}
		\No{S(t) \int_t^\I S(-\ta) Y(\ta) \ps(\ta) d\ta}_{L^4_t L^2_\om L^\mu_x} \\
		&\quad
		+
		\ph(\|V\|_{L^2_t \dH^{\si,c}_x})
		\|V\|_{L^2_t L^d_x}
		\No{ |\na|^\si S(t) \int_t^\I S(-\ta) Y(\ta) \ps(\ta) d\ta}_{L^4_t L^2_\om L^q_x} \\
		&\le
		\ph(\|V\|_{L^2_t \dH^{\si,c}_x})
		\|V\|_{L^2_t \dH^{\si,c}_x}
		\|h_f\|_{L^1(\BR)}
		\||\na|^{\si-\tw} \ps\|_{L^2_t L^2_x}.
	\end{align*}
	
	Next, we prove \eqref{TV+} when $\si \in [-\frac{d}{2}+\frac{3}{2}, \tw]$.
	Let $\ps(t,x) \in \CS(\BR^{1+d})$.
	Then we have
	\begin{align*}
		|\Br{\ps, |\na|^\si T_V(Z_0)}_{L^2_t L^2_x}|
		&= |\Br{|\na|^{\tw-\si} T_V^*(|\na|^\si \ps), |\na|^{\si-\tw}Z_0}_{L^2_\om L^2_x}| \\
		&\le \||\na|^{\tw-\si} T_V^*(|\na|^\si \ps)\|_{L^2_\om L^2_x}  \||\na|^{\si-\tw}Z_0\|_{L^2_\om L^2_x} \\
		&\le \ph(\|V\|_{L^2_t \dH^{1/2-\si,\wt{c}}_x}) \|\ps\|_\Ltx  \||\na|^{\si-\tw}Z_0\|_{L^2_\om L^2_x}.
	\end{align*}
Therefore, we have
\begin{align*}
	\||\na|^\si T_V(Z_0)\|_\Ltx \le \ph(\|V\|_{L^2_t \dH^{1/2-\si,\wt{c}}_x}) \||\na|^{\si-\tw}Z_0\|_{L^2_\om L^2_x}.
\end{align*}
Next, we prove \eqref{TV+} when $\si \in [\tw, \frac{d}{2}-\tw]$.
	\begin{align*}
		\||\na|^\si T_V(Z_0)\|_{L^2_t L^2_x}
		&= \||\na|^\si \BE[\overline{Y(t)} S_V(t)Z_0]\|_\Ltx \\
		&\le \||\na|^\si \BE[\overline{Y(t)} S(t)Z_0]\|_\Ltx \\
		&\quad + \No{|\na|^\si \BE \Dk{ \overline{Y(t)} \int_0^t S(t-\ta)(V(\ta)S_V(\ta)Z_0)d\ta} }_\Ltx \\
		&=: C + D.
	\end{align*}
	For $C$, we have by Theorem \ref{RS}
	\begin{align*}
		C &\lesssim \|h_f\|_{L^1(\BR)} \|Z_0\|_{L^2_\om \dot{H}^{\si-1/2}_x}.
	\end{align*}
Define $q$ and $\wt{\mu}$ by
	\begin{align}
	\frac{1}{q}=\frac{1}{2}-\frac{1}{2d},
	\quad \frac{1}{q'}=\frac{1}{\wt{c}} + \frac{1}{\wt{\mu}}.
\end{align}
By Proposition \ref{MA}, we have
	\begin{align*}
		D
		&\lesssim
		\|h_f\|_{L^1(\BR)}
		\|V(t)S_V(t)Z_0\|_{L^{4/3}_t L^2_\om \dot{H}^{\si-1/2,q'}_x} \\
	    &\lesssim
	    \|h_f\|_{L^1(\BR)}
	    \K{\||\na|^{\si-\tw} V\|_{L^2_t L^{\wt{c}}_x}
	    \|S_V(t)Z_0\|_{L^2_\om L^4_t L^{\wt{\mu}}_x}
	    +
	    \|V\|_{L^2_t L^d_x}
	    \||\na|^{\si-\tw} S_V(t)Z_0\|_{L^2_\om L^4_t L^q_x} } \\
    	&\le
    	\|h_f\|_{L^1(\BR)}
    	\ph(	\|V\|_{L^2_t \dH^{\si-1/2,\wt{c}}_x})
    	\|V\|_{L^2_t \dH^{\si-1/2,\wt{c}}_x}
    	\||\na|^{\si-\tw}Z_0\|_{L^2_\om L^2_x}.
	\end{align*}
Finally, \eqref{TV*+} when $\si \ge 0$ can be proved by the duality argument.
\end{proof}
We need a variants of \eqref{TV+}.
By the following proposition, we can replace the constant in \eqref{TV+} by $\ph(\|V\|_{L^2_t L^d_x})$.
\begin{prop.} \label{A4}
	Let $d \ge 3$ and $s = \frac{d-3}{2}$.
	Then there exist a monotone increasing function $\ph:[0,\I) \to [0,\I)$ such that 
	\begin{align*}
		\|T_V(Z_0)\|_{L^2_t \dH^{-s}_x}
		\le \ph(\|V\|_{L^2_t L^d_x}) \|h_f\|_{L^1(\BR)}\|Z_0\|_{L^2_\om(L^2_x \cap \dot{H}^{-s-1/2}_x)}.
	\end{align*}
\end{prop.}

\begin{proof}
	First, we have
	\begin{align*}
		&\|\BE[\ov{Y(t)} S_V(t)Z_0]\|_{L^2_t H^{-s}_x} \\
		&\quad \le \||\na|^{-s} \BE[\ov{Y(t)} S(t)Z_0]\|_{L^2_t L^2_x}
		+ \No{ |\na|^{-s} \BE\Dk{\ov{Y(t)}\int_0^t S(t-\ta)V(\ta)S_V(\ta)Z_0 d\ta }}_\Ltx \\
		&\quad =: A + B.
	\end{align*}
	For $A$, we have 
	\begin{align*}
		A \lesssim \|h_f\|_{L^1(\BR)}\||\na|^{-s-1/2} Z_0\|_{L^2_\om L^2_x}
	\end{align*}
	
	Next, we estimate $B$.
	Let $\frac{1}{r} = \frac{1}{2}-\frac{1}{d}.$
	By Proposition \ref{MA}, we have
	\begin{align*}
		&\No{ |\na|^{-s} \BE\Dk{\ov{Y(t)}\int_0^t S(t-\ta)V(\ta)S_V(\ta)Z_0 d\ta } }_\Ltx \\
		&\quad \lesssim \|h_f\|_{L^1(\BR)} \| |\na|^{-s-\tw} (V(t)S_V(t)Z_0) \|_{L^1_t L^2_x L^2_\om} \\
		&\quad \lesssim \|h_f\|_{L^1(\BR)} \|V(t) S_V(t)Z_0\|_{L^1_t L^2_\om L^{d'}_x} \\
		&\quad \le \|h_f\|_{L^1(\BR)} \|V\|_{L^2_t L^2_x} \|S_V(t)Z_0\|_{L^2_t L^2_\om L^r_x} \\
		&\quad \le \ph(\|V\|_{L^2_t L^d_x}) \|h_f\|_{L^1(\BR)} \|V\|_{L^2_t L^2_x}\|Z_0\|_{L^2_\om L^2_x}.
	\end{align*}
\end{proof}

Now we prove a useful lemma.
\begin{lem.}\label{AM}
	Let $d \ge 3$.
	Let $p,q \in [2,\I]$ satisfy $p >2$ and $\frac{2}{p}+\frac{d}{q} = \frac{d}{2}$.
	Let $\si$ satisfy $
		0 \le \si \le \frac{d}{2}-1
	$.
Let $\frac{1}{c} = \frac{\si+1}{d}$.
	Then there exists a monotone increasing function $\ph:[0,\I) \to [0,\I)$ such that the following holds:
	\begin{align*}
		\No{
			\int_0^t S_V(t,\ta)(\ps(\ta)Y(\ta)) d\ta
		}_{L^p_t L^2_\om \dot{H}^{\si,q}_x }
		\le
		\ph(\|V\|_{L^2_t \dH^{\si,c}_x})
		\|h_f\|_{L^1(\BR)}
		\|\ps\|_{\dot{H}^{\si-1/2}_x}.
	\end{align*}
\end{lem.}

\begin{proof}
	By Lemma \ref{SV} and Proposition \ref{fundamental lemma}, we have
	\begin{align*}
		&\No{
			\int_0^\I S_V(t,\ta)(\ps(\ta)Y(\ta)) d\ta
			}_{L^p_t L^2_\om \dot{H}^{\si,q}_x } \\
		&\quad = \No{|\na|^\si S_V(t)\int_0^\I S_V(\ta)^* (\ps(\ta)Y(\ta)) d\ta }_{L^2_\om L^p_t L^q_x} \\
		&\quad \le \ph(\|V\|_{L^2_t \dH^{\si,c}_x})
		         \No{|\na|^\si \int_0^\I S_V(\ta)^* (\ps(\ta)Y(\ta)) d\ta}_{L^2_\om L^2_x} \\
		&\quad \le \ph(\|V\|_{L^2_t \dH^{\si,c}_x})
		       \|h_f\|_{L^1(\BR)} \|\ps\|_{L^2_t \dot{H}^{\si-1/2}_x }.
	\end{align*}
Therefore, by the Christ-Kiselev's lemma, we have the desired result.
\end{proof}

\section{Nonlinear analysis}
\subsection{General case}
The estimates necessary to prove the Theorem \ref{rh} are collected in this subsection.
\begin{lem.}\label{A1}
		Let $d \ge 3$, $s = \frac{d-3}{2}$ and $\frac{1}{q} = \frac{1}{2} + \frac{s}{d}$.
		Then there exists a monotone increasing function $\ph:[0,\I) \to [0,\I)$ such that the following hold:
		\begin{align*}
			&\|S_V(t) Z_0\|_{L^2_\om L^4_t L^{2q}_x}
			\le
			\ph(\|V\|_{L^2_t(\dot{H}^{1/4}_x \cap \dH^{1/4,\I}_x)})
			\|Z_0\|_{L^2_\om \dot{H}^{1/4}_x}, \\
			& \|D_V(WY)\|_{L^4_t L^2_\om L^{2q}_x }
			\le
			\ph(\|V\|_{L^2_t(\dot{H}^{1/4}_x \cap \dH^{1/4,\I}_x)})
			\|h_f\|_{L^1(\BR)}
			\|W\|_{L^2_t \dot{H}^{-1/4}_x}.
		\end{align*}
\end{lem.}

\begin{proof}
	Let $\frac{1}{q_0} = \frac{1}{2} - \frac{1}{2d}$. 
	Note that 
	\begin{align*}
		\frac{2}{4}+\frac{d}{q_0} = \frac{d}{2},\quad  0 \le s \le \frac{d}{2}-1, 
		\quad \frac{1}{q_0}-\frac{s}{d} = \frac{1}{d}.
	\end{align*}
	By Lemma \ref{SV}, we have
	\begin{align*}
		\|S_V(t) Z_0\|_{L^2_\om L^4_t L^{2q}_x}
		&\lesssim
		\||\na|^\frac{1}{4} S_V(t)Z_0\|_{L^2_\om L^4_t L^{q_0}_x} \\
		&\le
	    \ph(\|V\|_{L^2_t(\dot{H}^{1/4}_x \cap \dH^{1/4,\I}_x)})
		\|Z_0\|_{L^2_\om \dot{H}^{1/4}_x}^2.
	\end{align*}

Next, we prove the second estimate.
We have
\begin{align*}
    \|D_V(WY)\|_{L^4_t L^2_\om L^{2q}_x }
	&\lesssim
	\No{|\na|^\frac{1}{4} \int_0^t S_V(t,\ta)(W(\ta)Y(\ta)) d\ta}_{L^4_t L^2_\om L^{q_0}_x}.
\end{align*}

By Lemma \ref{AM}, we have
\begin{align*}
	&\No{
		|\na|^{\ft} \int_0^t S_V(t,\ta)(W(\ta)Y(\ta)) d\ta
	}_{L^4_t L^2_\om L^{q_0}_x} \\
	&\quad
	\le
	\ph(\|V\|_{L^2_t(\dot{H}^{1/4}_x \cap \dH^{1/4,\I}_x)})
	\|h_f\|_{L^1(\BR)}
	\|W\|_{L^2_t \dot{H}^{-1/4}_x}.
\end{align*} 
\end{proof}

\begin{lem.} \label{A5}
	Let $d \ge 3$ and $s = \frac{d-3}{2}$.
	Then there exist a monotone increasing function $\ph:[0,\I) \to [0,\I)$ such that 
	\begin{align*}
		\|Q(V)\|_{L^2_t \dot{H}^{-s}_x}
		\le \ph(\|V\|_{L^2_t(L^2_x \cap L^\I_x)})
		\|h_f\|_{L^1(\BR)} \|V\|_{L^2_t L^2_x}  \|V\|_{L^2_t \dot{H}^{-1/2}_x}.
	\end{align*}
\end{lem.}

\begin{proof}
Let $\frac{1}{q}= \frac{1}{2}-\frac{1}{2d}$.
By Proposition \ref{MA}, we have
\begin{align*}
	&\No{ |\na|^{-s} \BE \Dk{ \ov{Y(t)} \int_0^t S(t-s) V(s) \int_0^s S_V(s,\ta) (V(\ta)Y(\ta))d\ta ds} }_\Ltx \\
	&\quad \lesssim \|h_f\|_{L^1(\BR)} \No{ V(t) \int_0^t S_V(t,\ta) (V(\ta)Y(\ta)) d\ta }_{L^{4/3}_t L^2_\om \dot{H}^{-s-1/2,q'}_x}.
\end{align*}
Then we have by Lemma \ref{AM}
\begin{align*}
	&\No{
		V(t)
		\int_0^t
		S_V(t,\ta)
		(V(\ta)Y(\ta))
		d\ta
	}_{L^{4/3}_t L^2_\om \dot{H}^{-s-1/2,q'}_x} \\
	&\quad
	\lesssim
	\No{
		V(t)
		\int_0^t
		S_V(t,\ta)
		(V(\ta)Y(\ta))
		d\ta
	}_{L^{4/3}_t L^2_\om L^{(2d)'}_x} \\
	&\quad
	\le
	\|V\|_{L^2_t L^2_x}
	\No{
		\int_0^t
		S_V(t,\ta)
		(V(\ta)Y(\ta))
		d\ta 
	}_{L^4_t L^2_\om L^q_x} \\
	&\quad
	\le
	\ph(\|V\|_{L^2_t(L^2_x \cap L^\I_x)})
	\|V\|_\Ltx
	\|V\|_{L^2_t \dot{H}^{-1/2}_x}.
\end{align*}
From the above, we have
\begin{align*}
	\|Q(V)\|_{L^2_t \dot{H}^{-s}_x}
	&\le
	\ph(\|V\|_{L^2_t(L^2_x \cap L^\I_x)})
	\|h_f\|_{L^1(\BR)}
	\|V\|_{\Ltx}
	\|V\|_{L^2_t \dot{H}^{-1/2}_x}.
\end{align*}
\end{proof}

\subsection{3D case}
In this subsection, we explain the strategy to prove Theorem \ref{rh 3D}
and collect the necessary estimates.
We use the cancellation used in \cite{2020CS}.
Thus, we transform the equation \eqref{rh2}
to more suitable form.
 
On the one hand, we have
\begin{align*}
	D_V(VY) &= D_0(VY) + D_0(V D_V(VY)).
\end{align*}
Thus, we have
\begin{align*}
	\BA_3(\rh) &= \BE[|D_V(VY)|^2] \\
	& = \BE[|D_0(VY)|^2] + 2\Re \BE[D_0(VY) \ov{D_0(V D_V(VY))}] + \BE[|D_0(V D_V(VY))|^2].
\end{align*}
We estimate the last two terms directly; however it seems difficult to estimate the first term itself.
On the other hand, we have
\begin{align*}
	\BA_5(\rh) &=Q(V) \\
	&= -2\Re \BE \Dk{ \ov{Y(t)} \int_0^t S(t-s) V(s) \int_0^s S(s-\ta) (V(\ta)Y(\ta))d\ta ds} \\
	&\quad -2\Re\BE \Dk{(-i)  \ov{Y(t)} \int_0^t S(t-s) V(s) \int_0^s S(s-r) V(r) \int_0^r S_V (r,\ta) (V(\ta)Y(\ta))d\ta  dr ds} \\
	&=:Q_2'(V) + Q_3(V).
\end{align*}
We can estimate $Q_3(V)$ directly.
Let $Q_2(V) = \BE[|D_0(VY)|^2] + Q_2'(V)$.
It seems difficult to estimate each of $\BE[|D_0(VY)|^2]$ and $Q_2'(V)$, 
but we can estimate $Q_2(V)$.

\begin{lem.}\label{A1 3D}
Let $d=3$.
Then there exists a monotone increasing function $\ph:[0,\I) \to [0,\I)$ such that the following holds:
\begin{align}
	&\|\BE[|S_V(t)Z_0|^2]\|_{L^2_t H^{1/2}_x}
	\le \ph(\|V\|_{L^2_t H^{1/2}_x}) \|Z_0\|_{L^2_t H^{1/2}_x}, \label{A1 3D a}\\
	&\|\BE[|S_{V_1}(t)Z_0|^2] - \BE[|S_{V_2}(t)Z_0|^2]\|_\Ltx  \nonumber \\
	&\qquad \le \ph(\|V_1\|_{L^2_t H^{1/2}_x} + \|V_2\|_{L^2_t H^{1/2}_x})
	\|Z_0\|_{L^2_t H^{1/2}_x}^2
	\|V_1-V_2\|_{L^2_t H^{1/2}_x} \label{A1 3D sa}.
\end{align}
\end{lem.}

\begin{proof}
First, we prove \eqref{A1 3D a}.
We have
\begin{align*}
	\|\BE[|S_V(t)Z_0|^2]\|_{L^2_t H^{1/2}_x}
	&\sim \|\BE[|S_V(t)Z_0|^2]\|_{L^2_t L^2_x} + \||\na|^\tw \BE[|S_V(t)Z_0|^2]\|_{L^2_t L^2_x} \\
	&=: A + B.
\end{align*}
For $A$, we have by Lemma \ref{SV}
\begin{align*}
	A
	&\le \|S_V(t)Z_0\|_{L^2_\om L^4_t L^4_x} \\
	&\le \|S_V(t)Z_0\|_{L^2_\om L^4_t H^{1/4,3}_x} \\
	&\le \ph(\|V\|_{L^2_t \dH^{1/4,12/5}_x}) \|Z_0\|_{L^2_\om \dH^{1/4}_x} \\
	&\le \ph(\|V\|_{L^2_t \dH^{1/2}_x}) \|Z_0\|_{L^2_\om \dH^{1/4}_x}.
\end{align*}
For $B$, we have by Lemma \ref{SV}
\begin{align*}
	B
	&\lesssim \||\na|^\tw S_V(t)Z_0\|_{L^2_\om L^4_t L^3_x}^2 \\
	&\le \ph(\|V\|_{L^2_t \dH^{1/2}_x}) \|Z_0\|_{L^2_\om \dH^{1/2}_x}.
\end{align*}

Next, we prove that there exists a monotone increasing function $\ph:[0,\I) \to [0,\I)$ such that the followings hold:
\begin{align}
	&\|S_{V_1}(t)Z_0 - S_{V_2}(t)Z_0\|_{L^2_\om L^4_t H^{1/2,3}_x} \nonumber \\
	&\quad \le \ph(\|V_1\|_{L^2_t H^{1/2}_x} + \|V_2\|_{L^2_t H^{1/2}_x})
	\|Z_0\|_{L^2_\om L^2_x}
	\|V_1-V_2\|_{L^2_t L^2_x}. \label{b sa 1}
\end{align}
Note that
\begin{align}\label{S_V t ta}
	S_{V_1}(t,\ta) - S_{V_2}(t,\ta) = -i\int_\ta^t S_{V_1}(t,r)(V_1(r)-V_2(r))S_{V_2}(r,\ta) dr.
\end{align}
\eqref{S_V t ta} and Lemma \ref{SV} imply
\begin{align*}
	&\|S_{V_1}(t)Z_0 - S_{V_2}(t)Z_0\|_{L^4_t L^2_\om L^3_x} \\
	&\quad = \No{ \int_0^t S_{V_1}(t,\ta) (V_1(\ta) - V_2(\ta)) S_{V_2}(\ta) Z_0 d\ta }_{L^4_t L^2_\om L^3_x} \\
	&\quad \le \ph(\|V_1\|_{L^2_t L^3_x})
	\No{(V_1(t) - V_2(t))S_{V_2}(t)Z_0 }_{L^{4/3}_t L^2_\om L^{3/2}_x} \\
	&\quad \le \ph(\|V_1\|_{L^2_t H^{1/2}_x} + \|V_2\|_{L^2_t H^{1/2}_x})
	\|Z_0\|_{L^2_\om L^2_x}
	\|V_1-V_2\|_{L^2_t L^2_x}.
\end{align*}
Moreover, we have by fractional Leibniz rule and Lemma \ref{SV}
\begin{align*}
	&\|S_{V_1}(t)Z_0 - S_{V_2}(t)Z_0\|_{L^4_t L^2_\om \dH^{1/2,3}_x} \\
	&\quad = \No{ |\na|^\tw \int_0^t S_{V_1}(t,\ta) (V_1(\ta) - V_2(\ta)) S_{V_2}(\ta) Z_0 d\ta }_{L^4_t L^2_\om L^3_x} \\
	&\quad \le \ph(\|V_1\|_{L^2_t H^{1/2}_x})
	\No{ |\na|^\tw (V_1(t) - V_2(t))S_{V_2}(t)Z_0 }_{L^{4/3}_t L^2_\om L^{3/2}_x} \\
	&\quad \le \ph(\|V_1\|_{L^2_t H^{1/2}_x} + \|V_2\|_{L^2_t H^{1/2}_x})
	\|Z_0\|_{L^2_\om H^{1/2}_x}
	\|V_1-V_2\|_{L^2_t H^{1/2}_x}.
\end{align*}
Therefore, we have by \eqref{b sa 1} and Lemma \ref{SV}
\begin{align*}
	&\|\BE[|S_{V_1}(t)Z_0|^2] - \BE[|S_{V_2}(t)Z_0|^2]\|_\Ltx \\
	&\quad \le \|S_{V_1}(t)Z_0 + S_{V_2}(t)Z_0\|_{L^4_t L^2_\om L^6_x}
	\|S_{V_1}(t)Z_0 - S_{V_2}(t)Z_0\|_{L^4_t L^2_\om L^3_x} \\
	&\quad \le \ph(\|V_1\|_{L^2_t H^{1/2}_x} + \|V_2\|_{L^2_t H^{1/2}_x})
	\|Z_0\|_{L^2_\om H^{1/2}_x}^2
	\|V_1-V_2\|_\Ltx.
\end{align*}
\end{proof}

\begin{lem.}\label{A2 3D}
	Let $d=3$.
	Then there exists a monotone increasing function $\ph:[0,\I) \to [0,\I)$ such that the following holds:
	\begin{align}
		&\|\BE[(S_V(t)Z_0)\overline{D_V(VY)}]\|_{\Ltx}
		\le \ph(\|V\|_{L^2_t \dH^{1/2}_x}) \|Z_0\|_{L^2_\om L^2_x} \|V\|_{L^2_t L^2_x}, \label{A2 3D a}\\
		&\|\BE[(S_{V_1}(t)Z_0)\ov{D_{V_1}(V_1 Y)}] - \BE[(S_{V_2}(t)Z_0)\ov{D_{V_2}(V_2 Y)}]\|_\Ltx \nonumber \\
		&\qquad \le \ph(\|V_1\|_{L^2_t H^{1/2}_x} + \|V_2\|_{L^2_t H^{1/2}_x})
		\|h_f\|_{L^1(\BR)}
		\|Z_0\|_{L^2_\om H^{1/2}_x}
		\|V_1-V_2\|_{L^2_t H^{1/2}_x}. \label{A2 3D sa}
	\end{align}
\end{lem.}

\begin{proof}
First, we prove \eqref{A2 3D a}.
	By Lemma \ref{SV} and Proposition \ref{fundamental lemma}, we have
	\begin{align*}
		\|\BE[(S_V(t)Z_0)\overline{D_V(VY)}]\|_{\Ltx}
		&\le
		\|S_V(t)Z_0\|_{L^2_\om L^4_t L^3_x}
		\|D_V(VY)\|_{L^4_t L^2_\om L^6_x} \\
		&\le
		\ph(\|V\|_{L^2_t L^d_x})
		\|Z_0\|_{L^2_\om L^2_x}
		\||\na|^\tw D_V(VY)\|_{L^4_t L^2_\om L^3_x} \\
		&\le \ph(\|V\|_{L^2_t \dH^{1/2}_x}) \|Z_0\|_{L^2_\om L^2_x} \|V\|_{L^2_t L^2_x}.
	\end{align*}

Next, we prove that there exists a monotone increasing function $\ph:[0,\I) \to [0,\I)$ such that the followings hold:
\begin{align}
	&\|D_{V_1}(V_1 Y) - D_{V_2}(V_2 Y)\|_{L^4_t L^2_\om \dH^{1/2,3}_x} \nonumber \\
	&\quad \le \ph(\|V_1\|_{L^2_t H^{1/2}_x} + \|V_2\|_{L^2_t H^{1/2}_x})
	\|h_f\|_{L^1(\BR)}
	(\|V_1\|_{L^2_t H^{1/2}_x} + \|V_2\|_{L^2_t H^{1/2}_x})
	\|V_1-V_2\|_{L^2_t H^{1/2}_x}. \label{b sa 2}
\end{align}
\eqref{S_V t ta} and Lemma \ref{SV} imply
\begin{align*}
	&\|D_{V_1}(V_1 Y) - D_{V_2}(V_2 Y)\|_{L^4_t L^2_\om \dH^{1/2,3}_x} \\
	&\quad \le
	\|D_{V_1}(V_1 Y) - D_{V_2}(V_1 Y)\|_{L^4_t L^2_\om H^{1/2,3}_x}
	+ \|D_{V_2}((V_1-V_2) Y)\|_{L^4_t L^2_\om H^{1/2,3}_x} \\
	&=: A_1 + A_2.
\end{align*}
For $A_1$, we have by Lemma \ref{SV}, the fractional Leibniz rule and Lemma \ref{AM}
\begin{align*}
	A_1
	&\le \ph(\|V_1\|_{L^2_t \dH^{1/2}_x})
	\No{|\na|^\tw (V_1(t)-V_2(t)) \int_0^t S_{V_2}(t,\ta)(V_1(\ta)Y(\ta))}_{L^{4/3}_t L^2_\om L^{3/2}_x} \\
	&\le \ph(\|V_1\|_{L^2_t \dH^{1/2}_x})
	\|V_1-V_2\|_{L^2_t \dH^{1/2}_x}
	\No{\int_0^t S_{V_2}(t,\ta)(V_1(\ta)Y(\ta)) d\ta}_{L^4_t L^2_\om H^{1/2,3}_x} \\
	&\le \ph(\|V_1\|_{L^2_t \dH^{1/2}_x} + \|V_2\|_{L^2_t \dH^{1/2}_x})
	\|V_1-V_2\|_{L^2_t \dH^{1/2}_x}.
\end{align*}
Lemma \ref{AM} implies
\begin{align*}
	A_2 \le
	\ph(\|V_1\|_{L^2_t \dH^{1/2}_x} + \|V_2\|_{L^2_t \dH^{1/2}_x})
	\|V_1-V_2\|_{L^2_t \dH^{1/2}_x}.
\end{align*}
Therefore, we have by \eqref{b sa 1}, \eqref{b sa 2} and Lemma \ref{SV}
\begin{align*}
	&\|\BE[(S_{V_1}(t)Z_0)\ov{D_{V_1}(V_1 Y)}] - \BE[(S_{V_2}(t)Z_0)\ov{D_{V_2}(V_2 Y)}]\|_\Ltx \\
	&\quad \le \|\BE[(S_{V_1}(t)Z_0 - S_{V_2}(t)Z_0)\ov{D_{V_1}(V_1 Y)}] \|_\Ltx \\
	&\qquad  + \|\BE[(S_{V_2}(t)Z_0)\ov{(D_{V_1}(V_1 Y) -  D_{V_2}(V_2 Y))}]\|_\Ltx \\
	&\quad \le \|S_{V_1}(t)Z_0 - S_{V_2}(t)Z_0\|_{L^4_t L^2_\om L^3_x}
	\|D_{V_1}(V_1 Y)\|_{L^4_t L^2_\om L^6_x} \\
	&\qquad  + \|S_{V_2}(t)Z_0\|_{L^4_t L^2_\om L^3_x}
	\|D_{V_1}(V_1 Y) -  D_{V_2}(V_2 Y)\|_{L^4_t L^2_\om L^6_x} \\
	&\quad \le \ph(\|V_1\|_{L^2_t H^{1/2}_x} + \|V_2\|_{L^2_t H^{1/2}_x})
	\|h_f\|_{L^1(\BR)}
	\|Z_0\|_{L^2_\om H^{1/2}_x}
	\|V_1-V_2\|_{L^2_t H^{1/2}_x}.
\end{align*}

\end{proof}

\begin{lem.} \label{A3 3D}
	Let $d=3$.
	Then there exists a monotone increasing function $\ph:[0,\I) \to [0,\I)$ such that the following hold:
\begin{align}
	&\|\BE[D_0(VY) \ov{D_0(V D_V(VY))}]\|_{L^2_t L^2_x}
	\le
	\ph(\|V\|_{L^2_t H^{1/2}_x})
	\|h_f\|_{L^1(\BR)}^2
	\|V\|_\Ltx^3, \label{A3 3D a1}\\
	&\|\BE[D_0(V_1Y) \ov{D_0(V_1 D_{V_1}(V_1Y))}] - \BE[D_0(V_2 Y) \ov{D_0(V_2 D_{V_2}(V_2Y))}] \|_\Ltx  \nonumber \\ 
	&\quad \le \ph(\|V_1\|_{L^2_t H^{1/2}_x} + \|V_2\|_{L^2_t H^{1/2}_x})
	\|h_f\|_{L^1(\BR)}^2
	(\|V_1\|_{L^2_t H^{1/2}_x}^2 + \|V_2\|_{L^2_t H^{1/2}_x}^2)
	\|V_1-V_2\|_{L^2_t H^{1/2}_x}, \label{A3 3D sa1}\\
	&\|\BE[|D_0(V D_V(VY))|^2]\|_{L^2_t L^2_x}
	\le
	\ph(\|V\|_{L^2_t H^{1/2}_x})
	\|h_f\|_{L^1(\BR)}^2
	\|V\|_\Ltx^2
	\|V\|_{L^2_t H^{1/4}_x}^2, \label{A3 3D a2} \\
	&\|\BE[|D_0(V_1 D_{V_1}(V_1 Y))|^2] - \BE[|D_0(V_2 D_{V_2}(V_2 Y))|^2]\|_\Ltx \nonumber  \\
	&\quad \le \ph(\|V_1\|_{L^2_t H^{1/2}_x} + \|V_2\|_{L^2_t H^{1/2}_x})
	\|h_f\|_{L^1(\BR)}^2
	(\|V_1\|_{L^2_t H^{1/2}_x}^3 + \|V_2\|_{L^2_t H^{1/2}_x}^3)
	\|V_1-V_2\|_{L^2_t H^{1/2}_x}. \label{A3 3D sa2}
\end{align}
\end{lem.}

\begin{proof}
We only prove \eqref{A3 3D a1} and \eqref{A3 3D a2}
because one can prove \eqref{A3 3D sa1} and \eqref{A3 3D sa2} in the same way.

First, we have by Lemma \ref{AM}
\begin{align*}
	&\|D_0(VY) \ov{D_0(V D_V(VY))}\|_{L^2_t L^1_\om L^2_x} \\ 
	&\quad
	\le
	\|D_0(VY)\|_{L^2_\om L^4_t L^6_x}
	\|D_0(V D_V(VY))\|_{L^4_t L^2_\om L^3_x} \\
	&\quad
	\le
	\||\na|^\tw D_0(VY)\|_{L^2_\om L^4_t L^3_x}
	\|D_0(V D_V(VY)\|_{L^4_t L^2_\om L^3_x} \\
	&\quad
	\lesssim
	\|h_f\|_{L^1(\BR)}
	\|V\|_{L^2_t L^2_x}
	\|V D_V(VY)\|_{L^{4/3}_t L^2_\om L^{3/2}_x} \\
	&\quad
	\lesssim
	\|h_f\|_{L^1(\BR)}
	\|V\|_\Ltx^2
	\|D_V(VY)\|_{L^4_t L^2_\om L^6_x} \\
	&\quad
	\lesssim
	\|h_f\|_{L^1(\BR)}
	\|V\|_\Ltx^2
	\||\na|^\tw D_V(VY)\|_{L^4_t L^2_\om L^3_x} \\
	&\quad
	\le
	\ph(\|V\|_{L^2_t \dH^{1/2}_x} )
	\|h_f\|_{L^1(\BR)}^2
	\|V\|_\Ltx^3.
\end{align*}

Next, we have also by Lemma \ref{AM}
\begin{align*}
	&\|D_0(V D_V(VY))\|_{L^4_t L^2_\om L^4_x} \\
	&\quad
	\lesssim
	\||\na|^\ft V D_V(VY)\|_{L^{4/3}_t L^2_\om L^{3/2}_x} \\
	&\quad
	\lesssim
	\||\na|^\ft V\|_\Ltx \|D_V(VY)\|_{L^4_t L^2_\om L^6_x}
	+
	\|V\|_{L^2_t L^{12/5}_x} \||\na|^\ft D_V(VY)\|_{L^4_t L^2_\om L^4_x} \\
	&\quad
	\lesssim
	\||\na|^\ft V\|_\Ltx \||\na|^\tw D_V(VY)\|_{L^4_t L^2_\om L^3_x}
	+
	\|V\|_{L^2_t L^{12/5}_x}
	\||\na|^\tw D_V(VY)\|_{L^4_t L^2_\om L^3_x} \\
	&\quad
	\le
	\ph(\|V\|_{L^2_t \dH^{1/2}_x})
	\|h_f\|_{L^1(\BR)}
	\|V\|_\Ltx
	\|V\|_{L^2_t \dH^{1/4}_x}.
\end{align*}
\end{proof}

\begin{lem.} \label{Q3 3D}
	Let $d=3$.
	Then there exists a monotone increasing function $\ph:[0,\I) \to [0,\I)$ such that the following holds:
	\begin{align}
		&\|Q_3(V)\|_\Ltx
		\le
		\ph(\|V\|_{L^2_t H^{1/2}_x})
		\|h_f\|_{L^1(\BR)}^2
		\|V\|_\Ltx^3, \label{Q3 3D a} \\
		&\|Q_3(V_1)-Q_3(V_2)\|_\Ltx  \nonumber \\
		&\quad \le \ph(\|V_1\|_{L^2_t H^{1/2}_x} + \|V_2\|_{L^2_t H^{1/2}_x})
		\|h_f\|_{L^1(\BR)}^2
		(\|V_1\|_{L^2_t H^{1/2}_x}^2 + \|V_2\|_{L^2_t H^{1/2}_x}^2)
		\|V_1-V_2\|_{L^2_t H^{1/2}_x}. \label{Q3 3D sa}
	\end{align}
\end{lem.}

\begin{proof}
We only prove \eqref{Q3 3D a}.
Define $F$ and $G$ by
\begin{align*}
&Q_3(V) =: -2\Re\BE \Dk{(-i)  \ov{Y(t)} \int_0^t S(t-s) V(s) F(s) ds} \\
&F(s) =: \int_0^s S(s-r) V(r) G(r) dr.
\end{align*}
By Propositions \ref{fundamental lemma}, \ref{MA}, 
\begin{align*}
	\|Q_3(V)\|_\Ltx
	&\le
	\No{ \BE \Dk{ \ov{Y(t)} \int_0^t S(t-s) (V(s)F(s)) ds } }_\Ltx \\
	&\lesssim
	\|h_f\|_{L^1(\BR)}
	\||\na|^{-\tw}V(t)F(t)\|_{L^1_t L^2_\om L^2_x} \\
	&\lesssim
	\|h_f\|_{L^1(\BR)}
	\|V\|_\Ltx
	\No{ \int_0^t S(t-s) (V(s) G(s))  ds }_{L^2_t L^2_\om L^6_x} \\
	&\lesssim
	\|h_f\|_{L^1(\BR)}
	\|V\|_\Ltx
	\No{V(t) \int_0^t S_V (t,\ta) (V(\ta)Y(\ta))d\ta }_{L^2_t L^2_\om L^{6/5}_x} \\
	&\lesssim
	\|h_f\|_{L^1(\BR)}
	\|V\|_\Ltx^2
	\No{\int_0^t S_V (t,\ta) (V(\ta)Y(\ta))d\ta }_{L^\I_t L^2_\om L^3_x}\\
	&\lesssim
	\|h_f\|_{L^1(\BR)}
	\|V\|_\Ltx^2
	\No{|\na|^\tw\int_0^t S_V (t,\ta) (V(\ta)Y(\ta))d\ta }_{L^\I_t L^2_\om L^2_x}.
\end{align*}
We have by Lemma \ref{SV}
\begin{align*}
	&\No{|\na|^\tw\int_0^\I S_V (t,\ta) (V(\ta)Y(\ta))d\ta }_{L^\I_t L^2_\om L^2_x} \\
	&\quad \le
	\ph(\|V\|_{L^2_t \dH^{1/2}_x})
	\No{|\na|^\tw \int_0^\I S_V (\ta)^* (V(\ta)Y(\ta))d\ta }_{L^2_\om L^2_x}\\
	&\quad \le
	\ph(\|V\|_{L^2_t \dH^{1/2}_x})
	\|h_f\|_{L^1(\BR)}
	\|V\|_\Ltx.
\end{align*}
Therefore, the Christ-Kiselev lemma implies
\begin{align}
	&\No{|\na|^\tw\int_0^t S_V (t,\ta) (V(\ta)Y(\ta))d\ta }_{L^\I_t L^2_\om L^2_x}
\le \ph(\|V\|_{L^2_t \dH^{1/2}_x})
	\|h_f\|_{L^1(\BR)}
	\|V\|_\Ltx,
\end{align}
which completes the proof.
\end{proof}

\begin{lem.} \label{Q2 3D}
	Let $d=3$. Then we have
	\begin{align*}
	\|Q_2(V)\|_\Ltx \lesssim A_\th(h_f) \|V\|_{L^2_t H^{1/2}_x}^2,
	\end{align*}
where $A_\th(h_f)$ is defined by \eqref{C1}
\end{lem.}

\begin{proof}
	The proof is essentially same as that of \cite[Lemma 5.6, Proposition 5.7]{2022CS}.
	Recall the notations of \cite{2022CS}:
	\begin{align}
		&K(\et,\et_2) = h_f(2t\et + 2s\et_2)\sin(t(|\et|^2-\et_2\cdot \et)) \sin(t\et_2 \cdot \et + s|\et_2|^2) \nonumber \\
		&Q_2(U,V):= 2\Re \BE \Dk{ \ov{W_V(Y)}W_U(Y) + \ov{Y}\K{ W_V(W_U(Y)) + W_U(W_V(Y))} }, \label{Q2}
	\end{align}
where $W_V(\phi)$ is $D(V\phi)$ in our notation.
By the same argument as the proof of \cite[Lemma 5.6]{2022CS}, we have
\begin{align*}
	\|K\|_{L^2_t L^1_s}^2
	    &\le
	        \frac{1}{|\et_2|^{1-2\th}\sqrt{|\et|^2|\et_2|^2-(\et\cdot \et_2)^2 }}
	        \int_\BR dv
	        \K{
	        	\int_\BR du
	        	|u|^\th \Abs{h_f(\sqrt{u^2+v^2})}
	        }^2.
\end{align*}
Therefore, by the same argument as the proof of \cite[Proposition 5.7]{2022CS}, we have
\begin{align*}
	\|\wh{Q_2}(U,V)\|_{L^2_t}
	&\le A_\th[h_f] \int_{\BR^3} d\et_2
	     |\et_2|^{\th-\tw}\K{ |\et|^2 |\et_2|^2 - (\et\cdot\et_2)^2}^{-\ft}\\
	&\quad \times \Dk{
	     	\|\wh{U}(\et-\et_2,\cdot)\|_{L^2_t} 
	     	\|\wh{V}(\et_2,\cdot)\|_{L^2_t}
	     	+
	     	\|\wh{V}(\et-\et_2,\cdot)\|_{L^2_t}
	     	\|\wh{U}(\et_2,\cdot)\|_{L^2_t}
     	}.
\end{align*}
By the above formula, our problem reduced to the estimate of the following integral:
\begin{align*}
I := \int_{\BR^3 \times \BR^3} d\et d\xi
	     \K{ |\et|^2 |\xi|^2 - (\et\cdot\xi)^2}^{-\ft}
	     u(\et-\xi) v(\xi) \phi(\et) \Br{\et-\xi}^{-\tw} \lxr^{-\tw}|\xi|^{\th-\tw},
\end{align*}
where $\phi$ is an arbitrary test function and
\begin{align*}
	u(\et):= \Br{\et}^\tw\|\wh{U}(\et,\cdot)\|_{L^2_t}, 
	\quad
	v(\et):= \Br{\et}^\tw\|\wh{V}(\et,\cdot)\|_{L^2_t}.
\end{align*}
Therefore, we have
\begin{align*}
	I &\le \int_{\BR^2 \times \BR^2} d\et' d\xi' |\xi_1 \et_2 - \xi_2 \et_1|^{-\tw} \\
	  &\quad \times \int_{\BR \times \BR} d\xi_3 d \xi_3
	      u(\et'-\xi',\et_3-\xi_3) v(\xi',\xi_3) \phi(\et',\et_3) \Br{\et-\xi}^{-\tw} \lxr^{-\tw}|\xi|^{\th-\tw} \\
	  &\le \int_{\BR^2 \times \BR^2} d\et' d\xi' |\xi_1 \et_2 - \xi_2 \et_1|^{-\tw}
	  \|u(\et'-\xi',\cdot)\|_{L^2_{\et_3}} \|v(\xi',\cdot)\|_{L^2_{\et_3}} \|\phi(\et',\cdot)\|_{L^2_{\et_3}}.
\end{align*}
The rest of the proof is completely same as that of \cite[Proposition 5.7]{2022CS}.
\end{proof}

\begin{lem.}\label{A5 3D}
	Let $d=3$.
	Then there exists a monotone increasing function $\ph:[0,\I) \to [0,\I)$ such that the following holds:
	\begin{align}
		&\|Q(V)\|_{L^2_t \dH^{1/2}_x}
		\le \ph(\|V\|_{L^2_t H^{1/2}_x}) \|V\|_\Ltx^2, \label{A5 3D a}\\
		&\|Q(V_1)-Q(V_2)\|_{L^2_t \dH^{1/2}_x} \nonumber \\
		&\quad \le 
		\ph(\|V_1\|_{L^2_t H^{1/2}_x} + \|V_2\|_{L^2_t H^{1/2}_x})
		\|h_f\|_{L^1(\BR)}^2
		\|V_1-V_2\|_{L^2_t H^{1/2}_x}. \label{A5 3D sa}
	\end{align}

\end{lem.}
\begin{proof}
	We only prove \eqref{A5 3D a}.
	We have by Propositions \ref{MA} and \ref{fundamental lemma}
	\begin{align*}
		\||\na|^\tw \BA_5(\rh)\|_\Ltx
		&\le \No{|\na|^\tw \BE \Dk{ \ov{Y(t)} \int_0^t S(t-s) V(s) \int_0^s S_V(s,\ta) (V(\ta)Y(\ta))d\ta ds} }_\Ltx \\
		&\lesssim \No{ V(t) \int_0^t S_V(t,\ta) (V(\ta)Y(\ta))d\ta }_{L^{4/3}_t L^2_\om L^{3/2}_x} \\
		&\lesssim \|V\|_\Ltx \No{ \int_0^t S_V(t,\ta) (V(\ta)Y(\ta))d\ta }_{L^4_t L^2_\om L^6_x} \\
		&\le \ph(\|V\|_{L^2_t \dH^{1/2}_x}) \|V\|_\Ltx^2.
	\end{align*}
\end{proof}

\section{Proof of the main results}
\begin{proof}[{\bf Proof of Theorem \ref{rh}}]
Define $\Ph$ and $E(R)$ by
\begin{align}
	&E(R):= \{ \rh \in L^2_t(\BR,\dH^{-s}_x) : \|\rh\|_{L^2_t\dH^{-s}_x} \le R\}, \\
	&\Ph[\rh] := (1-L)^{-1} \K{ \sum_{n=1}^5 \BA_n(\rh)},
\end{align}
where $R > 0$ will be chosen later. 
We prove $\Ph:E(R) \to E(R)$ is a contraction map.
In the following estimates, we interpret that constants depend on $f$ and $w$.

Let $\frac{1}{q} = \frac{1}{2}-\frac{1}{d}$ and  $\frac{1}{r} = \frac{1}{2} + \frac{s}{d}$.
Let $X := \dH^{s-\tw,1}_x \cap \dH^{s+\ft,1}_x \cap \dH^{s}_x \cap \dH^{s+\ft}_x$.

First, we prove $\Ph:E(R) \to E(R)$ is well-defined.
Let $\rh \in E(R)$.
Note that $(1-L)^{-1} \in \CB(\Ltx)$ is equivalent to $(1-L)^{-1} \in \CB(L^2_t H^{-s}_x)$ since $L$ is a Fourier multiplier.
Then we have
\begin{align*}
	\|\Ph[\rh]\|_{L^2_t \dH^{-s}_x}
	&\le C(f,w) \sum_{n=1}^5 \|\BA_n(\rh)\|_{L^2_t \dH^{-s}_x}.
\end{align*}
By Lemma \ref{A1}, we have
\begin{align*}
	\|\BA_1(\rh)\|_{L^2_t \dH^{-s}_x} 
	&= \|\BE[|S_V(t)Z_0|^2]\|_{L^2_t \dH^{-s}_x} \\
	&\lesssim \|\BE[|S_V(t)Z_0|^2]\|_{L^2_t L^r_x} \\
	&\le \ph(\|V\|_{L^2_t (\dot{H}^{1/4}_x \cap \dot{H}^{1/4,\I}_x)}) \|Z_0\|_{L^2_\om \dot{H}^{1/4}_x}^2 \\
	&\le \ph(\|w\|_{\dH^{s+1/4,1}_x \cap \dH^{s+1/4}_x} \|\rh\|_{L^2_t \dot{H}^{-s}_x}) \|Z_0\|_{L^2_\om \dot{H}^{1/4}_x}^2 \\
	&\le \ph(R\|w\|_X) \ep_0^2.
\end{align*}
By Lemma \ref{A1}, we have
\begin{align*}
	\|\BA_3(\rh)\|_{L^2_t \dH^{-s}_x}
	&= \|\BE[|D_V(VY)|^2]\|_{L^2_t \dH^{-s}_x} \\
	&\lesssim \|\BE[|D_V(VY)|^2]\|_{L^2_t L^r_x} \\
	&\le \ph(\|V\|_{L^2_t(\H)}) \|V\|_{L^2_t \dot{H}^{-1/4}_x}^2 \\
	&\le \ph(\|w\|_{\dot{H}^{s+1/4,1}_x \cap \dot{H}^{s+1/4}_x} \|\rh\|_{L^2_t \dH^{-s}_x })
         \|w\|_{\dot{H}^{s-1/4,1}_x}^2 \|\rh\|_{L^2_t \dot{H}^{-s}_x}^2 \\
    &\le \ph(R\|w\|_X)R^2\|w\|_X^2 .
\end{align*}
By the Cauchy-Schwarz inequality, we have
\begin{align*}
	\|\BA_2(\rh)\|_{L^2_t \dH^{-s}_x}
	&\lesssim \|\BA_2(\rh)\|_{L^2_t L^r_x} \\
	&\le \|\BA_1(\rh)\|_{L^2_t L^r_x}^{\frac{1}{2}} \|\BA_3(\rh)\|_{L^2_t L^r_x}^{\frac{1}{2}}.
\end{align*}
By Proposition \ref{A4}, we have
\begin{align*}
	\|\BA_4(\rh)\|_{L^2_t \dH^{-s}_x}
	&\le \|\BE[\ov{Y(t)}S_V(t)Z_0]\|_{L^2_t \dH^{-s}_x} \\
	&\le
	\ph(
	    \|V\|_{L^2_t (L^2_x \cap L^\I_x)}
	    )
	\|Z_0\|_{L^2_\om(L^2_x \cap \dot{H}^{-s-1/2}_x)} \\
    &\le
    \ph(
        \|w\|_{\dot{H}^{s,1}_x \cap \dot{H}^{s}_x}
        \|\rh\|_{L^2_t \dot{H}^{-s}_x}
        )
    \|Z_0\|_{L^2_\om(L^2_x \cap \dot{H}^{-s-1/2}_x)} \\
    &\le \ph(R\|w\|_X) \ep_0.
\end{align*}
By Lemma \ref{A5}, we have
\begin{align*}
	\|\BA_5(\rh)\|_{L^2_t \dH^{-s}_x} 
	&\le \|Q(V)\|_{L^2_t L^2_x} \\
	&\le  \ph(\|V\|_{L^2_t(L^2_x \cap L^\I_x)})
	      \|V\|_\Ltx \|V\|_{L^2_t \dot{H}^{-1/2}_x}	\\
    &\le \ph(\|w\|_{\dH^{s,1}_x\cap \dH^{s}_x} \|\rh\|_{L^2_t\dH^{-s}_x})
         \|w\|_{\dH^{s,1}_x} \|w\|_{\dH^{s-1/2,1}_x} \|\rh\|_{L^2_t \dH^{-s}_x}^2 \\
    &\le \ph(R\|w\|_X) R^2 \|w\|_X^2 .
\end{align*}
From the above, if we choose $R>0$ and $\ep_0 > 0$ sufficiently small, we have
\begin{align*}
	\|\Ph[\rh]\|_{L^2_t \dH^{-s}_x} \le R.
\end{align*}

Now we prove the contraction.
Let $\rh_1, \rh_2 \in E(R)$ and $V_1:= w \ast \rh_1$, $V_2:= w \ast \rh_2$.
We have
\begin{align*}
	&\|\Ph[\rh_1]-\Ph[\rh_2]\|_{L^2_t \dH^{-s}_x} \\
	&\quad \le
	     \No{
	     	\BE\Dk{|S_{V_1}(t)Z_0+D_{V_1}(V_1 Y)|^2
	     	-
	        |S_{V_2}(t)Z_0 + D_{V_2}(V_2 Y)|^2}
     	}_{L^2_t \dH^{-s}_x} \\ 
	&\qquad +
	     \|\BE[\ov{Y(t)}(S_{V_1}(t)-S_{V_2}(t))Z_0]\|_{L^2_t \dH^{-s}_x}
	     +
	     \|Q(V_1)-Q(V_2)\|_{L^2_t \dH^{-s}_x}\\
	&\quad \le
	     \|S_{V_1}(t)Z_0+D_{V_1}(V_1 Y) + S_{V_2}(t)Z_0+D_{V_2}(V_2 Y)\|_{L^4_t L^2_\om L^{2r}_x} \\
	&\qquad \times \|S_{V_1}(t)Z_0+D_{V_1}(V_1 Y)-S_{V_2}(t)Z_0-D_{V_2}(V_2 Y)\|_{L^4_t L^2_\om L^{2r}_x} \\
	&\qquad + \|\BE[\ov{Y(t)}(S_{V_1}(t)-S_{V_2}(t))Z_0]\|_{L^2_t \dH^{-s}_x} + \|Q(V_1)-Q(V_2)\|_{L^2_t \dH^{-s}_x} \\
	&\quad \le
	    CR\K{\|S_{V_1}(t)Z_0-S_{V_2}(t)Z_0\|_{L^4_t L^2_\om L^{2r}_x}
	    	+ \|D_{V_1}(V_1 Y)-D_{V_2}(V_2 Y)\|_{L^4_t L^2_\om L^{2r}_x}
	       } \\
    &\qquad +\|\BE[\ov{Y(t)}(S_{V_1}(t)-S_{V_2}(t))Z_0]\|_{L^2_t \dH^{-s}_x} + \|Q(V_1)-Q(V_2)\|_{L^2_t \dH^{-s}_x}
\end{align*} 
By \eqref{S_V t ta}, we have
\begin{align*}
	&\|S_{V_1}(t)Z_0-S_{V_2}(t)Z_0\|_{L^4_t L^2_\om L^{2r}_x} \nonumber \\
	&\quad \le \ph(\|V_1\|_{L^2_t H^{1/4,3}_x}+\|V_2\|_{L^2_t H^{1/4,3}_x})
	\|V_1-V_2\|_{L^2_t H^{1/4,3}_x}.
\end{align*}
And also, we have
\begin{align*}
	&\|D_{V_1}(V_1 Y)-D_{V_2}(V_2 Y)\|_{L^4_t L^2_\om L^{2r}_x} \\
	&\quad \le
	    \No{ \int_0^t(S_{V_1}(t,\ta)-S_{V_2}(t,\ta))(V_1(\ta)Y(\ta))d\ta }_{L^4_t L^2_\om L^{2r}_x} \\
	&\qquad  + \No{ \int_0^t S_{V_1}(t,\ta)(V_1(\ta)-V_2(\ta))Y(\ta)d\ta}_{L^4_t L^2_\om L^{2r}_x} \\
	&\quad =:A_1 + A_2.
\end{align*}
The problem is the estimate of $A_1$.
By \eqref{S_V t ta}, we have
\begin{align*}
	A_1
	&\le \No{
	    	\int_0^t S_{V_1}(t,r)(V_1(r)-V_2(r))
	    	\int_0^r S_{V_2}(r,\ta) (V_1(\ta)Y(\ta)) d\ta dr
    	}_{L^4_t L^2_\om L^{2r}_x} \\
    &\le \|V_1-V_2\|_{L^2_t (H^{1/4}_x \cap H^{1/4,\I}_x)}
         \No{|\na|^\tw \int_0^t S_{V_2}(t,\ta)(V_1(\ta)Y(\ta)) d\ta}_{L^2_t L^2_\om L^q_x}.
\end{align*}
We can estimate $\|\BE[\ov{Y(t)}(S_{V_1}(t)-S_{V_2}(t))Z_0]\|_{L^2_t \dH^{-s}_x}$ and $\|Q(V_1)-Q(V_2)\|_{L^2_t \dH^{-s}_x}$ in the same way as the above.
The uniqueness of the solution can be shown in the same way.
\end{proof}

\begin{proof}[{\bf Proof of Theorem \ref{rh 3D}}]
Define $\Ph$ and $E(R)$ by
\begin{align}
	&E(R):= \{ \rh \in L^2_t(\BR, H^\tw_x) : \|\rh\|_{L^2_t H^{1/2}_x} \le R\}, \\
	&\Ph[\rh] := (1-L)^{-1} \K{ \sum_{n=1}^5 \BA_n(\rh)},
\end{align}
where $R > 0$ will be chosen later. 
We prove $\Ph:E(R) \to E(R)$ is a contraction map.

In the following estimates, we interpret that constants depend on $f$ and $w$.

First, we estimate $\|\Ph[\rh]\|_\Ltx$.
By Lemma \ref{A1 3D}, we have
\begin{align*}
	\|\BA_1(\rh)\|_\Ltx 
	&\le
	\ph(\|\rh\|_{L^2_t H^{1/2}_x})
	\|Z_0\|_{L^2_\om H^{1/2}_x}^2 
	\le \ph(R) \ep_0^2.
\end{align*}
By Lemma \ref{A2 3D}, we have
\begin{align*}
	\|\BA_2(\rh)\|_\Ltx
	&\le
	\ph(\|\rh\|_{L^2_t H^{1/2}_x})
	\|Z_0\|_{L^2_\om L^2_x}
	\|\rh\|_{L^2_t L^2_x}
	\le \ph(R) R \ep_0.
\end{align*}
By Proposition \ref{A4}, we have
\begin{align*}
	\|\BA_4(\rh)\|_{L^2_t L^2_x}
	&\le \ph(\|\rh\|_{L^2_t H^{1/2}_x}) \|Z_0\|_{L^2_\om(L^2_x \cap \dot{H}^{-1/2}_x)}
	\le \ph(R) \ep_0.
\end{align*}
Now we have
\begin{align*}
	\BA_3(\rh) + \BA_5(\rh)
	=
	2\Re \BE[D_0(VY) \ov{D_0(V D_V(VY))}] + \BE[|D_0(V D_V(VY))|^2]+ Q_2(V)+Q_3(V).
\end{align*}
By Lemma \ref{A3 3D}, we have
\begin{align*}
	&\|2\Re \BE[D_0(VY) \ov{D_0(V D_V(VY))}]\|_\Ltx
	\le
	\ph(\|V\|_{L^2_t H^{1/2}_x})
	\|V\|_\Ltx^3 
	\le
	\ph(R)
	R^3,\\
	&\|\BE[|D_0(V D_V(VY))|^2]\|_{L^2_t L^2_x}
	\le
	\ph(\|V\|_{L^2_t H^{1/2}_x })
	\|V\|_{L^2_t H^{1/2}_x}^4
	\le 
	\ph(R) R^4.
\end{align*}
By Lemma \ref{Q3 3D},
we have
\begin{align*}
	\|Q_3(V)\|_\Ltx
	&\le
	\ph(\|\rh\|_{L^2_t H^{1/2}_x})
	\|\rh\|_\Ltx^3
	\le
	\ph(R) R^3.
\end{align*}
By Lemma \ref{Q2 3D}, we have
\begin{align*}
	\|Q_2(V)\|_\Ltx &\lesssim \|V\|_{L^2_t H^{1/2}_x}^2
	\lesssim R^2.
\end{align*}

Next, we estimate $\||\na|^\tw \Ph[\rh]\|_\Ltx$.
By Lemma \ref{A1 3D}, we have
\begin{align*}
	\||\na|^\tw \BA_1(\rh)\|_\Ltx
	\le \ph(\|V\|_{L^2_t \dH^{1/2}_x}) \|Z_0\|_{L^2_t H^{1/2}_x}
	\le \ph(R)\ep_0.
\end{align*}
By the fractional Leibniz rule and Proposition \ref{fundamental lemma}, we have
\begin{align*}
	\||\na|^\tw \BA_3(\rh)\|_\Ltx
	&\lesssim \||\na|^\tw D_V(VY)\|_{L^4_t L^2_\om L^3_x}^2 \\
	&\le \ph(\|V\|_{L^2_t H^{1/2}_x}) \|V\|_\Ltx^2 \le \ph(R) R^2.
\end{align*}
We can estimate $\||\na|^\tw \BA_2(\rh)\|_\Ltx$ by the Cauchy-Schwarz inequality.
By Proposition \ref{fundamental lemma}, we have
\begin{align*}
	\||\na|^\tw \BA_4(\rh)\|_\Ltx \le \ph(\|V\|_{L^2_t L^3_x}) \|Z_0\|_{L^2_\om L^2_x}
	\le \ph(R) \ep_0.
\end{align*}
Finally, we have by Lemma \ref{A5 3D}
\begin{align*}
	\||\na|^\tw \BA_5(\rh)\|_\Ltx
	\le \ph(\|V\|_{L^2_t H^{1/2}_x}) \|V\|_\Ltx^2 
	\le \ph(R) R^2.
\end{align*}

From the above, if we choose $R>0$ and $\ep_0 > 0$ sufficiently small, we have
\begin{align*}
	\|\Ph[\rh]\|_{L^2_t H^{1/2}_x} \le R.
\end{align*}
Thus, $\Ph:E(R) \to E(R)$ is well-defined.

Next, we prove that $\Ph$ is a contraction map.
Let $\rh_1, \rh_2 \in E(R)$ and $V_1 := w \ast \rh_1$, $V_2:= w \ast \rh_2$.
First, we estimate $\|\Ph[\rh_1] - \Ph[\rh_2]\|_\Ltx$.
We have
\begin{align*}
	\|\Ph[\rh_1] - \Ph[\rh_2]\|_\Ltx
	&\le
	    \|\BA_1[\rh_1]-\BA_1[\rh_2]\|_\Ltx
	    + \|\BA_2[\rh_1]-\BA_2[\rh_2]\|_\Ltx \\
	&\quad +\|\BA_4[\rh_1]-\BA_4[\rh_2]\|_\Ltx+\|\BA_3[\rh_1]+\BA_5[\rh_1]-\BA_3[\rh_2]-\BA_5[\rh_2]\|_\Ltx.    
\end{align*}
For the first term, we have by Lemma \ref{A1 3D}
\begin{align*}
	&\|\BA_1[\rh_1]-\BA_1[\rh_2]\|_\Ltx
	= \|\BE[|S_{V_1}(t)Z_0|^2] - \BE[|S_{V_2}(t)Z_0|^2]\|_\Ltx \nonumber \\
	&\quad \le \ph(\|V_1\|_{L^2_t H^{1/2}_x} + \|V_2\|_{L^2_t H^{1/2}_x})
	\|Z_0\|_{L^2_t H^{1/2}_x}^2
	\|V_1-V_2\|_{L^2_t H^{1/2}_x} \\
	&\quad \le \ph(R) R^2 \|\rh_1-\rh_2\|_{L^2_t H^{1/2}_x}.
\end{align*}
For the second term, we have by Lemma \ref{A2 3D}
\begin{align*}
	&\|\BA_2[\rh_1] - \BA_2[\rh_2]\|_\Ltx
	\le \|\BE[(S_{V_1}(t)Z_0-S_{V_2}(t)Z_0)\ov{D_{V_1}(V_1 Y)}]\|_\Ltx \\
	&\qquad \le \ph(\|V_1\|_{L^2_t H^{1/2}_x} + \|V_2\|_{L^2_t H^{1/2}_x})
	\|h_f\|_{L^1(\BR)}
	\|Z_0\|_{L^2_\om H^{1/2}_x}
	\|V_1-V_2\|_{L^2_t H^{1/2}_x} \\
	&\qquad \le \ph(R)R\|\rh_1-\rh_2\|_{L^2_t H^{1/2}_x} 
\end{align*}
For the third term, it is easy to see that 
\begin{align*}
	&\|\BA_4[\rh_1]-\BA_4[\rh_2]\|_\Ltx \\
	&\quad \le \ph(\|V_1\|_{L^2_t H^{1/2}_x} + \|V_1\|_{L^2_t H^{1/2}_x})
	\|h_f\|_{L^1(\BR)}
	\|Z_0\|_{L^2_\om H^{1/2}_x}
	\|V_1 - V_2\|_{L^2_t H^{1/2}_x} \\
	&\quad \le \ph(R) R \|V_1 - V_2\|_{L^2_t H^{1/2}_x}.
\end{align*}
Finally, we have
\begin{align*}
	&\|\BA_3[\rh_1]+\BA_5[\rh_1]-\BA_3[\rh_2]-\BA_5[\rh_2]\|_\Ltx \\
	&\quad \le \|Q_2(V_1) - Q_2(V_2)\|_\Ltx + \|Q_3(V_1)-Q_3(V_2)\|_\Ltx \\
	&\qquad + \|\BE[D_0(V_1Y) \ov{D_0(V_1 D_{V_1}(V_1Y))}] - \BE[D_0(V_2 Y) \ov{D_0(V_2 D_{V_2}(V_2Y))}] \|_\Ltx  \\ 
    &\qquad + \|\BE[|D_0(V_1 D_{V_1}(V_1 Y))|^2] - \BE[|D_0(V_2 D_{V_2}(V_2 Y))|^2]\|_\Ltx.
\end{align*}
For the first term, we have
\begin{align*}
	\|Q_2(V_1) - Q_2(V_2)\|_\Ltx
	&\le \tw \|Q_2(V_1, V_1) - Q_2(V_2,V_2)\|_\Ltx \\
	&\le \|Q_2(V_1, V_1-V_2)\|_\Ltx + \|Q_2(V_1-V_2,V_2)\|_\Ltx,
\end{align*}
where $Q_2(\cdot, \cdot)$ was defined by \eqref{Q2}.
In the proof of Lemma \ref{Q2 3D}, we have proved
\begin{align*}
	\|Q_2(U,V)\|_\Ltx \lesssim \|U\|_{L^2_t H^{1/2}_x} \|V\|_{L^2_t H^{1/2}_x}.
\end{align*}
Thus, we have
\begin{align*}
	\|Q_2(V_1)-Q(V_2)\|_\Ltx
	&\lesssim \K{\|V_1\|_{L^2_t H^{1/2}_x} + \|V_2\|_{L^2_t H^{1/2}_x}}
	          \|V_1-V_2\|_{L^2_t H^{1/2}_x}.
\end{align*}
The last three terms can be estimated by Lemmas \ref{Q3 3D} and \ref{A3 3D}.
Therefore, we have
\begin{align*}
	\|\Ph[\rh_1] - \Ph[\rh_2]\|_\Ltx
	&\le \ph(R)(R+\ep_0)\|\rh_1-\rh_2\|_{L^2_t H^{1/2}_x} \\
	&\le \frac{1}{4} \|\rh_1-\rh_2\|_{L^2_t H^{1/2}_x},
\end{align*}
for sufficiently small $R, \ep_0 >0$.

Next, we estimate $\|\Ph[\rh_1] - \Ph[\rh_2]\|_{L^2_t \dH^{1/2}_x}$.
\begin{align*}
	&\|\Ph[\rh_1]-\Ph[\rh_2]\|_{L^2_t \dH^{1/2}_x} \\
	&\quad \le
	\No{
		\BE\Dk{|S_{V_1}(t)Z_0+D_{V_1}(V_1 Y)|^2
			-
			|S_{V_2}(t)Z_0 + D_{V_2}(V_2 Y)|^2}
	}_{L^2_t \dH^{1/2}_x} \\ 
	&\qquad +
	\|\BE[\ov{Y(t)}(S_{V_1}(t)-S_{V_2}(t))Z_0]\|_{L^2_t \dH^{1/2}_x}
	+
	\|Q(V_1)-Q(V_2)\|_{L^2_t \dH^{1/2}_x}\\
	&\quad =: A + B + C.
\end{align*}
For $A$, we have
\begin{align*}
	A
	&\le \No{ |\na|^\tw \K{S_{V_1}(t)Z_0+D_{V_1}(V_1 Y) + S_{V_2}(t)Z_0 + D_{V_2}(V_2 Y)}}_{L^4_t L^2_\om L^3_x}\\
	&\quad \times 
	\No{|\na|^\tw \K{S_{V_1}(t)Z_0+D_{V_1}(V_1 Y) - S_{V_2}(t)Z_0 - D_{V_2}(V_2 Y)}}_{L^4_t L^2_\om L^3_x} \\
	&\le \ph(R)(R+\ep_0)
	\K{
	       \|S_{V_1}(t)Z_0 - S_{V_2}(t)Z_0\|_{L^4_t L^2_\om \dH^{1/2,3}_x}
	       + \|D_{V_1}(V_1 Y) - D_{V_2}(V_2 Y)\|_{L^4_t L^2_\om \dH^{1/2,3}_x}
	 } \\
 &=  \ph(R)(R+\ep_0)(A_1 + A_2).
\end{align*}
We have by \eqref{b sa 1}
\begin{align*}
	A_1 \le \ph(\|V_1\|_{L^2_t H^{1/2}_x} + \|V_2\|_{L^2_t H^{1/2}_x})
	\|Z_0\|_{L^2_\om L^2_x}
	\|V_1-V_2\|_{L^2_t L^2_x}.
\end{align*}
And also, we have by \eqref{b sa 2}
\begin{align*}
	A_2 \le \ph(\|V_1\|_{L^2_t H^{1/2}_x} + \|V_2\|_{L^2_t H^{1/2}_x})
	\|h_f\|_{L^1(\BR)}
	(\|V_1\|_{L^2_t H^{1/2}_x} + \|V_2\|_{L^2_t H^{1/2}_x})
	\|V_1-V_2\|_{L^2_t H^{1/2}_x}.
\end{align*}
By \eqref{S_V t ta} and Proposition \ref{fundamental lemma}, we have 
\begin{align*}
	B \le \ph(\|V_1\|_{L^2_t H^{1/2}_x} + \|V_2\|_{L^2_t H^{1/2}_x}) \|Z_0\|_{L^2_\om H^{1/2}_x} \|V_1-V_2\|_{L^2_t H^{1/2}_x}.
\end{align*}
We have by Lemma \ref{A5 3D}
\begin{align*}
	C \le \ph(\|V_1\|_{L^2_t H^{1/2}_x} + \|V_2\|_{L^2_t H^{1/2}_x})
	\|h_f\|_{L^1(\BR)}^2
	\|V_1-V_2\|_{L^2_t H^{1/2}_x}. 
\end{align*}
From the above, we have
\begin{align*}
	\|\Ph[\rh_1] - \Ph[\rh_2]\|_{L^2_t \dH^{1/2}_x}
	&\le \ph(R)(R+\ep_0)(A_1 + A_2) + B + C \\
	&\le \ph(R)(R+\ep_0) \|\rh_1 - \rh_2\|_{L^2_t H^{1/2}_x} \\
	&\le \frac{1}{4} \|\rh_1 - \rh_2\|_{L^2_t H^{1/2}_x}
\end{align*}
for sufficiently small $R, \ep_0 >0$.

From the above, we have
\begin{align*}
	\|\Ph[\rh_1] - \Ph[\rh_2]\|_{L^2_t H^{1/2}_x}
	&\le \frac{1}{2} \|\rh_1 - \rh_2\|_{L^2_t H^{1/2}_x}
\end{align*}
for sufficiently small $R, \ep_0 > 0$.

The uniqueness of the solution can be shown in the same way.
\end{proof}

\section{Proof of the main lemmas}
In this section, we give a proof of Lemmas \ref{scattaring lemma} and \ref{scattaring lemma 1.5}.
In the following, we interpret constants depend on $f$.
First, we prove the following preliminary lemma.

\begin{lem.}\label{Lemma 1}
	Let $d \ge 3$.
	Let $0 \le \si \le \frac{d}{2}-1$ and $\frac{1}{c}=\frac{\si+1}{d}$.
	Let $Z_0 \in L^2_\om \dH^\si_x$.
	Assume that
	\begin{align*}
		V \in L^2_t(\BR, \dH^{\si-1/2}_x \cap \dH^{\si,c}_x).
	\end{align*}
	Let $h_f \in L^1(\BR)$.
	Define $Z(t)$ by $Z(t) := S_V(t)Z_0 + D_V(VY)(t)$.
	Then we have $Z(t) \in C(\BR,L^2_\om \dH^\si_x)$.
\end{lem.}

\begin{proof}
	First, we prove that 
	\begin{align} \label{SV conti}
		\BR \ni t \mapsto S_V(t)Z_0 \in L^2_\om \dH^\si_x \mbox{ is continuous. }
	\end{align}
It is clear that $\BR \ni t \mapsto S(t)Z_0 \in L^2_\om \dH^\si_x$ is continuous.
Therefore, the problem is the Duhamel term.
Let $t_n \to t_\I$ as $n \to \I$.
Then we have
\begin{align*}
	&\No{\int_0^{t_\I}
		S(t_\I-\ta)(V(\ta)S_V(\ta)Z_0)d\ta - \int_0^{t_n} S(t_n-\ta)(V(\ta)S_V(\ta)Z_0) d\ta}_{L^2_\om \dH^\si_x} \\
	&\quad \le \No{ (S(t_\I)-S(t_n))\int_0^{t_\I} S(\ta)^*|\na|^\si (V(\ta)S_V(\ta)Z_0) d\ta }_{L^2_\om L^2_x} \\
	&\qquad + \No{ \int_{t_n}^{t_\I} S(t_n-\ta) |\na|^\si (V(\ta)S_V(\ta)Z_0)d \ta }_{L^2_\om L^2_x} \\
	&\quad =: A_n + B_n. 
\end{align*}
By the fractional Leibniz rule and Lemma \ref{SV}, we have
\begin{align*}
	\int_0^{t_\I} \| S(\ta)^*|\na|^\si (V(\ta)S_V(\ta)Z_0)\|_{L^2_\om L^2_x} d\ta < \I.
\end{align*}
Thus, we have $A_n \to 0$ as $n \to \I$.
Similarly, $B_n \to 0$ as $n \to \I$.

Next, we consider the Duhamel term $D_V(VY)(t)$.
Let $t_n \to t_\I$ as $n \to \I$.
	We have
	\begin{align*}
		&\No{
			\int_0^{t_\I} S_V(t_\I,\ta)(V(\ta)Y(\ta))d\ta
			-
			\int_0^{t_n} S_V(t_n,\ta)(V(\ta)Y(\ta))d\ta
		}_{L^2_\om \dH^\si_x} \\
		&\quad
		\le
		\No{
			(S_V(t_\I)-S_V(t_n))\int_0^{t_\I} S_V(\ta)^*(V(\ta)Y(\ta)) d\ta
		}_{L^2_\om \dH^\si_x}
		+
		\No{
			\int_{t_n}^{t_\I} S_V(t_n,\ta)(V(\ta)Y(\ta)) d\ta
		}_{L^2_\om \dH^\si_x} \\
		&\quad =: C_n + D_n.
	\end{align*}
	By \eqref{TV*+}, we have
	\begin{align*}
		\No{\int_0^{t_\I} S_V(\ta)^*(V(\ta)Y(\ta)) d\ta}_{L^2_\om \dH^\si_x}
		&\le
		\ph(\|V\|_{L^2_t \dH^{\si,c}_x})
		\|V\|_{L^2_t \dH^{\si-1/2}_x} \\
		&< \I.
	\end{align*}
Thus, we have by \eqref{SV conti} that $C_n \to 0$ as $n \to \I$.
Since
\begin{align*}
	S_V(t_n,\ta) = S(t_n-\ta) - i\int_\ta^{t_n} S(t_n-r)V(r)S_V(r,\ta) dr,
\end{align*}
we have
\begin{align*}
	D_n
	&\le
	\No{
		\int_{t_n}^{t_\I} S(t_n-\ta)(V(\ta)Y(\ta)) d\ta
	}_{L^2_\om \dH^\si_x} \\
&\quad +
\No{
\int_{t_n}^{t_\I} \int_\ta^{t_n} S(t_n-r)V(r)S_V(r,\ta) dr (V(\ta)Y(\ta)) d\ta
}_{L^2_\om \dH^\si_x} \\
&=:D_n^1 + D_n^2.
\end{align*}
By \eqref{TV*+}, we have
\begin{align*}
	D_n^1 &\lesssim \|V\|_{L^2_t([t_n,t_\I],\dH^{\si-1/2}_x)} \\
	&\to 0 \mbox{  as  } n \to \I.
\end{align*}
Let $\frac{1}{q}:= \frac{1}{2}-\frac{1}{d}$ and $\frac{1}{p}:=\tw-\frac{\si}{d}$.
For $D_n^2$, we have by the fractional Leibniz rule and Lemma \ref{AM}
\begin{align*}
	D_n^2
	&=
	\No{
		\int_{t_\I}^{t_n} S(r)^* V(r) \int_{t_\I}^r S_V(r,\ta)(V(\ta)Y(\ta)) d\ta dr 
	}_{L^2_\om \dH^\si_x} \\
&\le
\No{
	|\na|^\si \Dk{ V(t) \int_{t_\I}^t S_V(t,\ta)(V(\ta)Y(\ta)) d\ta }
}_{L^2_\om L^2_t([t_\I,t_n], L^{q'}_x)} \\
&\le
\||\na|^\si V\|_{L^2_t([t_\I,t_n], L^c_x)}
\No{
	\int_{t_\I}^t S_V(t,\ta)(V(\ta)Y(\ta)) d\ta
}_{L^2_\om L^\I_t L^p_x} \\
&\qquad + \|V\|_{L^2_t([t_n,t_\I], L^d_x)} 
         \No{
         	|\na|^\si \int_0^t S_V(t,\ta) (V(\ta)Y(\ta)) d\ta
         }_{L^2_\om L^\I_t L^2_x} \\
&\to 0 \mbox{ as } n \to \I.
\end{align*}
\end{proof}

\begin{proof}[{\bf Proof of Lemma \ref{scattaring lemma}}]
Thanks to Lemma \ref{Lemma 1}, it suffices to prove the scattering.
We have
\begin{align*}
	&\|S(-t)Z(t)-S(-s)Z(s)\|_{L^2_\om H^\si_x} \\
	&\quad \le \|S(-t)S_V(t)Z_0 - S(-s)S_V(s)Z_0\|_{L^2_\om H^\si_x} \\
	&\qquad + \|S(-t)D_V(VY)(t) - S(-s)D_V(VY)(s)\|_{L^2_\om H^\si_x}.
\end{align*}

Let $\frac{1}{q}:= \frac{1}{2}-\frac{1}{d}$.
For the first term, we have by the endpoint Strichartz estimates
\begin{align}
	&\|S(-t)S_V(t)Z_0 - S(-s)S_V(s)Z_0\|_{L^2_\om H^\si_x} \nonumber \\
	&\quad
	=
	\No{ \int_s^t S(-\ta)V(\ta)S_V(\ta)Z_0 d\ta }_{L^2_\om H^\si_x} \nonumber \\
	&\quad
	\lesssim
	\|\sd^\si V(\ta)S_V(\ta)Z_0\|_{L^2_\ta([s,t], L^2_\om L^{q'}_x ) } \nonumber \\
	&\quad
	\lesssim
	\|V\|_{L^2_\ta([s,t],H^{\si}_x \cap H^{\si,\I}_x)}
	\|S_V(\ta)Z_0\|_{ L^\I_\ta L^2_\om H^\si_x} \nonumber \\
	&\quad
	\le
	\ph(\|V\|_{L^2_t (H^\si_x \cap H^{\si,\I}_x)})
	\|V\|_{L^2_\ta([s,t],H^{\si}_x \cap H^{\si,\I}_x)}
	\|Z_0\|_{L^2_\om H^\si_x} \label{yokutsukau} \\
	&\quad \to 0 \mbox{ as } t,s \to \pm \I. \nonumber 
\end{align}

For the second term, we have
\begin{align*}
	&\|S(-t)D_V(YV)(t)-S(-s)D_V(YV)(s)\|_{L^2_\om H^{\si}_x} \\
	&\quad \le \No{\int_s^t S(-t)S_V(t-\ta)(Y(\ta)V(\ta)) d\ta}_{L^2_\om H^{\si}_x} \\
	&\qquad + \No{\int_0^s (S(-t)S_V(t) - S(-s)S_V(s))S_V(\ta)^* (Y(\ta)V(\ta)) d\ta}_{L^2_\om H^{\si}_x} \\
	&\quad =:A(t,s) + B(t,s).
\end{align*}
By Lemma \ref{SV} and Proposition \ref{fundamental lemma}, we have
\begin{align*}
	A(t,s)
	&\le \ph(\|V\|_{L^2_t (H^\si_x \cap H^{\si,\I}_x)}) \No{\int_s^t S_V(\ta)^* (Y(\ta)V(\ta)) d\ta}_{L^2_\om H^\si_x} \\
	&\le \ph(\|V\|_{L^2_t (H^\si_x \cap H^{\si,\I}_x)})
	\|V\|_{L^2_t([s,t],\dH^{-1/2}_x \cap \dH^{\si-1/2}_x)} \\
	&\to 0 \;\; \mbox{as} \;\; t,s \to \pm \I.
\end{align*}
By the estimate for the first term, we have
\begin{align*}
	\|S(-t)S_V(t)- S(-s)S_V(s)\|_{\CB(L^2_\om H^\si_x)} \to 0 \;\; \mbox{as} \;\; t,s \to \pm \I.
\end{align*}
Therefore, we have by Proposition \ref{fundamental lemma}
\begin{align*}
	B(t,s) &\le \No{(S(-t)S_V(t) - S(-s)S_V(s)) \Dk{ \int_0^s S_V(\ta)^* (Y(\ta)V(\ta)) d\ta}}_{L^2_\om H^{\si}_x} \\
	&\le
	\|S(-t)S_V(t) - S(-s)S_V(s)\|_{\CB(L^2_\om H^\si_x)}
	\No{\int_0^s S_V(\ta)^* (Y(\ta)V(\ta)) d\ta}_{L^2_\om H^\si_x} \\
	&\le
	\|S(-t)S_V(t) - S(-s)S_V(s)\|_{\CB(L^2_\om H^\si_x)}
	\ph(\|V\|_{L^2_t (H^\si_x \cap H^{\si,\I}_x)})
	\|V\|_{L^2_t(\dH^{-1/2}_x \cap \dH^{\si-1/2}_x)}  \\
	&\to 0 \;\; \mbox{as} \;\; t,s \to \pm \I.
\end{align*}
\end{proof}

\begin{proof}[{\bf Proof of Lemma \ref{scattaring lemma 1.5}}]
Thanks to Lemma \ref{Lemma 1}, it suffices to prove the scattering.
We have
\begin{align*}
	&\|S(-t)Z(t) - S(-s)Z(s)\|_{L^2_\om \dH^{1/2}_x} \\
	&\quad
	\le
	 \|S(-t)S_V(t)Z_0-S(-s)S_V(s)Z_0\|_{L^2_\om \dH^{1/2}_x} \\
	&\qquad
	+
	 \|S(-t)D_V(VY)(t)-S(-s)D_V(VY)(s)\|_{L^2_\om \dH^{1/2}_x} \\
	&\quad =:C + D.
\end{align*}
For $C$, we have by the same estimate that led to \eqref{yokutsukau} 
\begin{align}
	C &\le 
	\|V\|_{L^2([s,t],\dH^{1/2}_x)}
	\ph(\|V\|_{L^2_t \dH^{1/2}_x})
	\|Z_0\|_{L^2_\om \dH^{1/2}_x} \nonumber\\
	&\to 0 \mbox{ as } t,s \to \I.  \label{norm bound} 
\end{align}
For $D$, we have
\begin{align*}
	&\|S(-t)D_V(VY)(t)-S(-s)D_V(VY)(s)\|_{L^2_\om \dH^{1/2}_x} \\
	&\quad
	\le
	\No{\int_s^t S(-t)S_V(t,\ta)(V(\ta)Y(\ta)) d\ta}_{L^2_\om \dH^{1/2}_x} \\
	&\qquad
	+
	\No{\int_0^s (S(-t)S_V(t) - S(-s)S_V(s))S_V(\ta)^* (V(\ta)Y(\ta)) d\ta}_{L^2_\om \dH^{1/2}_x} \\
	&\quad
	=:D_1 + D_2.
\end{align*}
For $D_1$, we have by \eqref{TV*+}
\begin{align*}
	D_1
	&= \No{\int_s^t S_V(\ta)^* (V(\ta)Y(\ta)) d\ta}_{L^2_\om \dH^{1/2}_x} \\
	&\le \|h_f\|_{L^1(\BR)} \ph(\|V\|_{L^2_t \dH^{1/2}_x}) \|V\|_{L^2_t([s,t],L^2_x)} \\
	&\to 0 \mbox{ as } t,s \to \I.
\end{align*}
Note that by \eqref{norm bound}, we have
\begin{align*}
	\|S(-t)S_V(t)-S(-s)S_V(s)\|_{\CB(L^2_\om \dH^{1/2}_x)} \to 0 \mbox{ as } t,s \to \I.
\end{align*}
Therefore, we have by \eqref{TV*+}
\begin{align*}
	D_2
	&\le
	\|S(-t)S_V(t)-S(-s)S_V(s)\|_{\CB(L^2_\om \dH^{1/2}_x)} 
	\No{\int_0^s S_V(\ta)^* (Y(\ta)V(\ta)) d\ta }_{L^2_\om \dH^{1/2}_x} \\
	&\le
	\|S(-t)S_V(t)-S(-s)S_V(s)\|_{\CB(L^2_\om \dH^{1/2}_x)} 
	\|h_f\|_{L^1(\BR)}
    \ph(\|V\|_{L^2_t \dH^{1/2}_x})
	\|V\|_\Ltx \\
	&\to 0 \mbox{ as } t,s \to \I.
\end{align*}
\end{proof}

\end{document}